%% file: tagawa.tex
\documentclass[12pt]{article}
\usepackage{makeidx}
\usepackage{latexsym}
\usepackage{amsfonts}
\usepackage{amsmath}
\usepackage{amsthm}
\usepackage{amssymb}
\usepackage{amscd}
\usepackage[dvips]{graphicx}
\usepackage{color}
\usepackage{eepic}
\usepackage{setspace}
\usepackage{calrsfs}
\theoremstyle{definition}
\newtheorem{theorem}{Theorem}[section]
\newtheorem{prop}[theorem]{Proposition}
\newtheorem{lemma}[theorem]{Lemma}
\newtheorem{corollary}[theorem]{Corollary}
\newtheorem{definition}[theorem]{Definition}

\newtheorem{remark}[theorem]{Remark}
\newtheorem{conjecture}[theorem]{Conjecture}
\newenvironment{demo}[1]{%
  \trivlist
  \item[\hskip\labelsep
        {\bf #1.}]
}{%
\hfill\qedsymbol
  \endtrivlist
}
 \makeatletter
    
    \@addtoreset{equation}{section}
  \makeatother
\renewcommand{\mathcal}{\mathrsfs}
\title{
Pfaffian decomposition
and a Pfaffian analogue of $q$-Catalan Hankel determinants
}
\author{
Masao ISHIKAWA\footnote{
Department of Mathematics,
Faculty of Education, Tottori University, Koyama, Tottori 680-8550, Japan,
{\tt ishikawa@fed.tottori-u.ac.jp}
}, \ 
Hiroyuki TAGAWA\footnote{
Faculty of Education, Wakayama University, Sakaedani, Wakayama 640-8510, Japan,
{\tt tagawa@math.edu.wakayama-u.ac.jp}
} \ 
and 
Jiang ZENG\footnote{
Institut Camille Jordan, Universit\'e Claude Bernard Lyon 1, 69622 Villeurbanne cedex, France,
{\tt zeng@math.univ-lyon1.fr}
}
}
\date{
}
\def\defterm#1{{\sl #1}\/}

\newcommand\Pf{\operatorname{Pf}}

\newcommand\sgn{\operatorname{sgn}}

\def\qbinom#1#2{\left[{{#1} \atop {#2}}\right]_{q}}
\newcommand\rdots{\mathinner{\mkern1mu\raise0pt\vbox{\kern7pt\hbox{.}}
     \mkern2mu\raise4pt\hbox{.}\mkern2mu\raise8pt\hbox{.}\mkern1mu}}
\def\V#1#2#3#4{V^{{#1},{#2}}\left({#3};{#4}\right)}

\def\covered{\mathinner{\mkern1mu\raise0pt\vbox{\kern7pt\hbox{$<$}}
     \mkern-4mu\raise2pt\hbox{.}\mkern2mu}}
\def\covers{\mathinner{\mkern1mu\raise0pt\vbox{\kern7pt\hbox{$>$}}
     \mkern-12mu\raise2pt\hbox{.}\mkern8mu}}

\def\PATH#1#2{{\cal P}\left({#1},{#2}\right)}
\def\NPATH#1#2{{\cal P}_0\left({#1},{#2}\right)}
\def\GF#1{{\operatorname{GF}}\left[{#1}\right]}

\input roby4.tex

\begin{document}
%
%
%
%% --------- < *** >---------- %%
%% Title
%% --------- < *** >---------- %%
\maketitle
%
%
%
%% --------- < *** >---------- %%
%% Abstract
%% --------- < *** >---------- %%
\begin{abstract}
Motivated by the Hankel determinant evaluation of moment sequences, we  study a kind of Pfaffian analogue evaluation.
We prove an $LU$-decomposition analogue for skew-symmetric matrices, called Pfaffian decomposition.
We then apply this formula to evaluate Pfaffians related to some moment sequences of classical orthogonal polynomials. 
%% 09/11/2010 %% In particular we evaluate  a kind of 
%% 09/11/2010 %% $q$-Catalan Pfaffians.
%% 09/11/2010 %% We also establish an application to weighted enumeration of shifted reverse plane partitions.
In particular we obtain a product formula for a kind of q-Catalan Pfaffians.  We also establish a connection between our Pfaffian formulas and 
certain  weighted enumeration of shifted reverse plane partitions.
\end{abstract}
%
%
%
%% --------- < *** >---------- %%
%% Introduction
%% --------- < *** >---------- %%
\input tagawa01.tex

%
%
%% --------- < *** >---------- %%
%% Pfaffian decomposition
%% --------- < *** >---------- %%
\input tagawa02.tex

%
%
%% --------- < *** >---------- %%
%% An Pfaffain analogue of $q$-Catalan Hankel determinants
%% --------- < *** >---------- %%
\input tagawa03.tex
%
%
%% --------- < *** >---------- %%
%% A proof
%% --------- < *** >---------- %%
\input tagawa04.tex

%
%
%% --------- < *** >---------- %%
%% An application
%% --------- < *** >---------- %%
\input tagawa05.tex
%
%
%% --------- < *** >---------- %%
%% Open problems
%% --------- < *** >---------- %%
\input tagawa06.tex

%
%
%% --------- < *** >---------- %%
%% Appendix
%% --------- < *** >---------- %%
\input tagawa07.tex
\end{document}

%% file: roby4.tex
\newdimen\Squaresize \Squaresize=20pt
\newdimen\thickness \thickness=1pt         
                                                    
\def\Square#1{\hbox{\vrule width \thickness
   \vbox to \Squaresize{\hrule height \thickness\vss                            
      \hbox to \Squaresize{\hss#1\hss}
   \vss\hrule height\thickness} 
\unskip\vrule width \thickness} 
\kern-\thickness}                                                            
                               
\def\vsquare#1{\vbox{\Square{$#1$}}\kern-\thickness}
\def\blank{\omit\hskip\Squaresize}

\def\fibyoung#1{\let\\=\cr              %added  let\\ =\cr 
\vbox{\smallskip\offinterlineskip
\halign{&\vsquare{##}\cr #1}}\,}

\def\borderlessrect#1#2{\hbox{%\hskip \thickness
   \vbox to \Squaresize{\vskip \thickness \vss
      \hbox to #2 {\hss #1\hss}
   \vss\vskip\thickness} 
\unskip \hskip -\thickness
} 
\kern-\thickness}                                                            
                               
\def\vborderlessrect#1#2{\vbox{\borderlessrect{$#1$}{#2}}\kern-\thickness}

\def\borderless#1{\omit\vborderlessrect{#1}{\Squaresize}}
\def\borderlessrc#1#2{\omit\vborderlessrect{#1}{#2}}

\def\msquare#1{\vbox{\hbox{\vrule width \thickness
   \vbox to \Squaresize{\hrule height \thickness
      \hbox to \Squaresize{\hfil{\sevenrm #1}}
   \vfil\hrule height\thickness}
\unskip\vrule width \thickness}
\kern-\thickness}\kern-\thickness}

\def\twosquare#1#2{\vbox{\hbox{\vrule width \thickness %\vskip 1pt
   \vbox to \Squaresize{\hrule height \thickness
      \hbox to \Squaresize{\hfil{\rm #1}}\vss
      \hbox to \Squaresize{\hss{#2}\hss}
   \vfil\hrule height\thickness}
\unskip\vrule width \thickness}
\kern-\thickness}\kern-\thickness}

\def\twoblank#1#2{\vbox{\hbox{
   \vbox to \Squaresize{\vskip 2pt
      \hbox to \Squaresize{\hfil{\sevenrm #1}\ }\vss
      \hbox to \Squaresize{\hss{#2}\hss}
   \vfil}\unskip\kern-\thickness}
}\unskip\kern-\thickness}

\def\young#1{
\def\>{\blank}
\def\<{\borderless}
\def\*{\borderlessrc}
\def\p{\omit\msquare}
\def\t{\omit\twosquare}
\def\b{\omit\twoblank}
\let\\=\cr %added  let\\ =\cr 
\vbox{\smallskip\offinterlineskip
\halign{&\vsquare{##}\cr #1}}}

\newdimen\smsquaresize \smsquaresize=12pt
\newdimen\smthickness \smthickness=.5pt
\font\smcellfont=cmss8 scaled \magstep0

\def\smsquare#1{\hbox{\vrule width \smthickness
   \unskip\vbox to \smsquaresize{\hrule height \smthickness\vss
      \hbox to \smsquaresize{\hss{\smcellfont #1}\hss}
   \vss\hrule height\smthickness} 
\unskip\vrule width \smthickness} 
\kern-\smthickness}

\def\smvsquare#1{\vbox{\smsquare{$#1$}}\kern-\smthickness}
\def\blank{\omit\hskip\smsquaresize}

\def\smyoung#1{\let\\=\cr %added  let\\ =\cr
\vbox{\smallskip\offinterlineskip
\halign{&\smvsquare{##}\cr #1}}}
\newdimen\vsmsquaresize \vsmsquaresize=10pt
\newdimen\vsmthickness \vsmthickness=.5pt
\font\vsmcellfont=cmsl8 scaled \magstep0
\font\vsmletterfont=cmr6 scaled \magstep0

\def\vsmsquare#1{\hbox{\vrule width \vsmthickness
   \unskip\vbox to \vsmsquaresize{\hrule height \vsmthickness\vss
      \hbox to \vsmsquaresize{\hss{\vsmcellfont #1}\hss}
   \vss\hrule height\vsmthickness} 
\unskip\vrule width \vsmthickness} 
\kern-\vsmthickness}
\def\vsmvsquare#1{\vbox{\vsmsquare{#1}}\kern-\vsmthickness}
\def\vsmblank{\omit\hskip\vsmsquaresize}
\def\vsmborderless#1{\hbox{\hskip \vsmthickness\unskip
   \vbox to \vsmsquaresize{\vss
      \hbox to \vsmsquaresize{\hss{\vsmletterfont #1}\hss}
   \vss} 
\unskip\hskip \vsmthickness} 
\kern-\vsmthickness}                                                            \def\vsmvborderless#1{\vbox{\vsmborderless{#1}}\kern-\vsmthickness}

\def\vsmyoung#1{
\def\>{\vsmblank}
\def\<{\omit\vsmvborderless}
\let\\=\cr %added  let\\ =\cr
\vbox{\smallskip\offinterlineskip
\halign{&\vsmvsquare{##}\cr #1}}}

\newdimen\edgesize \edgesize=20pt
\newdimen\eedgesize \eedgesize=21pt
\newdimen\doublesize \doublesize=41pt
\newdimen\ddoublesize \ddoublesize=43pt

\newdimen\triplesize \triplesize=62pt
\newdimen\tetrasize \tetrasize=85pt
\newdimen\futosa \futosa=1pt         

\def\Hako#1{\hbox{\vrule width \futosa
   \vbox to \eedgesize{\hrule height \futosa\vss                            
      \hbox to \edgesize{\hss#1\hss}
   \vss\hrule height\futosa} 
\unskip\vrule width \futosa} 
\kern-\futosa}                                                            
                               
\def\Seihokei#1{\vbox{\Hako{#1}}\kern-\futosa}

\def\Horizontal#1{\hbox{\vrule width \futosa
   \vbox to \eedgesize{\hrule height \futosa\vss                            
      \hbox to \doublesize{\hss#1\hss}
   \vss\hrule height \futosa} 
\unskip\vrule width \futosa} 
\kern-\futosa}                                                            
                               
\def\Vertical#1{\hbox{\vrule width \futosa
   \vbox to \ddoublesize{\hrule height \futosa\vss                            
      \hbox to \edgesize{\hss#1\hss}
   \vss\hrule height \futosa} 
\unskip\vrule width \futosa} 
\kern-\futosa}                                                            
                               
\def\NS#1{\vbox{\Hako{$#1$}\vskip22pt}\kern-\futosa}
\def\NSS#1{\vbox{\Hako{$#1$}\vskip42pt}\kern-\futosa}
\def\NSSS#1{\vbox{\Hako{$#1$}\vskip64pt}\kern-\futosa}
\def\H#1{\vbox{\Horizontal{$#1$}}\kern-\futosa}
\def\HH#1#2{\vbox{\Horizontal{$#1$}\Horizontal{$#2$}}\kern-\futosa}
\def\VER#1{\vbox{\Vertical{$#1$}}\kern-\futosa}
\def\VS#1{\vbox{\Vertical{$#1$}\vskip20pt}\kern-\futosa}
\def\VSS#1{\vbox{\Vertical{$#1$}\vskip42pt}\kern-\futosa}
\def\VV#1#2{\vbox{\Vertical{$#1$}\Vertical{$#2$}}\kern-\futosa}
\def\NV#1#2{\vbox{\Hako{$#1$}\Vertical{$#2$}}\kern-\futosa}
\def\NVS#1#2{\vbox{\Hako{$#1$}\Vertical{$#2$}\vskip22pt}\kern-\futosa}
\def\YTT#1#2#3{\vbox{\Horizontal{$#1$}\hbox{\vbox{\Vertical{$#2$}}\kern-\futosa\vbox{\Vertical{$#3$}}\kern-\futosa}}\kern-\futosa}

\def\domino#1{
\def\ns{\omit\NS}
\def\nss{\omit\NSS}
\def\nsss{\omit\NSSS}
\def\h{\omit\H}
\def\hh{\omit\HH}
\def\V{\omit\VER}
\def\Vs{\omit\VS}
\def\Vss{\omit\VSS}
\def\vsss{\omit\VSSS}
\def\vv{\omit\VV}
\def\nv{\omit\NV}
\def\nvs{\omit\NVS}
\def\ytt{\omit\YTT}
\let\\=\cr %added  let\\ =\cr 
\vbox{\smallskip\offinterlineskip
\halign{&\Seihokei{##}\cr #1}}}

%% file: tagawa01.tex
%%%% ---------------< ooo >--------------- %%%%
%%%%
%%%% Section 1
%%%%
%%%% ---------------< ooo >--------------- %%%%

\section{Introduction}

The Hankel determinants 
%% 01/11/2010 %% $
%% 01/11/2010 %% \det\Bigl(C_{i+j-2}\Bigr)_{1\leq i,j\leq n}
%% 01/11/2010 %% $
of Catalan numbers
have drawn the interests of many researchers with relations to 
the combinatorial arguments of lattice paths in recent years 
(see, e.g., \cite{A,Ba,BCQY,GX,ITZ,LT,Sta}).
%% it seems hard to find the same arguments on Pfaffians of Catalan numbers.
It is well-known that 
if $\{\mu_{n}\}_{n\geq0}$ is the moment sequence of certain orthogonal polynomials,
the Hankel determinant
$
\det\Bigl(\mu_{i+j-2}\Bigr)_{1\leq i,j\leq n}
$
have a nice formula
because of the classical theory of orthogonal polynomials
(see \cite{ITZ}).
In this paper we would like to exploit a Pfaffian analogue of this kind of Hankel determinants.

\smallskip
%
%
%% ------------------------- %%
%% skew-symmetric matix
%% ------------------------- %%
%
%
We say a matrix $A=(a^{i}_{j})_{i,j\geq1}$ 
(resp. $A=(a^{i}_{j})_{1\leq i,j\leq n}$) 
is \defterm{skew-symmetric} 
if  it satisfies $a^{j}_{i}=-a^{i}_{j}$ for $i,j\geq1$ (resp. $1\leq i,j\leq n$).
%% If we are given an $n\times n$ skew-symmetric matrix $A=(a^{i}_{j})_{1\leq i,j\leq n}$
%% of even degree,
%% then $\det A$ is the square of a polynomial of its entries $a^{i}_{j}$.
%% So the \defterm{Pfaffian} of $A$,
%% denoted by $\Pf A$,
%% is defined to be the square root of $\det A$,
%% where we take the branch which takes the value $\Pf J_{2n}=1$.
%% If $n$ be an even integer,
%% the \defterm{pfaffian} of $A$ is defined by
If we are given an $n\times n$ skew-symmetric matrix $A=(a^{i}_{j})_{1\leq i,j\leq n}$
where $n$ is an even integer,
then the \defterm{Pfaffian} of $A$ (see \cite{Sta,Ste}),
denoted by $\Pf A$,
is defined to be
\begin{equation}
\Pf(A)=\sum \epsilon(\sigma_{1},\sigma_{2},\hdots,\sigma_{n-1},\sigma_{n})\,
a^{\sigma_{1}}_{\sigma_{2}} \dots a^{\sigma_{n-1}}_{\sigma_{n}},
\label{eq:Pfaffian}
\end{equation}
where the summation is over all partitions $\sigma=\{\{\sigma_{1},\sigma_{2}\}_{<},\hdots,\{\sigma_{n-1},\sigma_{n}\}_{<}\}$
of $[n]:=\{1,2,\dots,n\}$ into $2$-elements blocks,
and where $\epsilon(\sigma_{1},\sigma_{2},\hdots,\sigma_{n-1},\sigma_{n})$ denotes the sign of the permutation
\begin{equation}
\begin{pmatrix}
1&2&\cdots&n-1&n\\
\sigma_{1}&\sigma_{2}&\cdots&\sigma_{n-1}&\sigma_{n}
\end{pmatrix},
\label{eq:1-factor}
\end{equation}
and we use the notation $[n]$.
A partition $\sigma$ of $[n]$ into $2$-elements blocks is called a \defterm{perfect} matching or a \defterm{$1$-factor}.
\smallskip

%% We insist on this view point so that entries of Pfaffians should be related to
%% moment functions of certain orthogonal polynomials,
%% and we look for the right definition of ``Hankel Pfaffians''.
%% One candidate is the Pfaffians defined in Pfaffian transforms,
%% i.e.,
%% $
%% \Pf\begin{pmatrix}
%% \mu_{j-i-1}
%% \end{pmatrix}_{1\leq i<j\leq 2n}
%% $.
%% Here we propose 
%% Here we propose another candidate,
%% i.e.
%% $
%% \Pf\Bigl((j-i)\mu_{i+j-2}\Bigr)_{1\leq i,j\leq 2n}
%% $
%% as a Pfaffian analogue of the above Hankel determinants,
%% and for a $q$-analogue we take
%% $
%% \Pf\Bigl((q^{i-1}-q^{j-1})\mu_{i+j-2}\Bigr)_{1\leq i,j\leq 2n}
%% $.
As most of the orthogonal polynomials have their $q$-analogues, 
in order to propose a  Pfaffian analogue of the above Hankel determinants of moments,
 we have  the ordinary  version 
and $q$-version.  More precisely, we propose $\Pf\bigl((j-i)\,\mu_{i+j+r-2}\bigr)_{1\leq i,j\leq2n}$ as a Pfaffian
analogue of the above Hankel determinants 
of the  moments $\mu_n$, and for a $q$-analogue of $\mu_n(q)$ of $\mu_n$ we take 
\[
\Pf\bigl((q^{i-1}-q^{j-1})\, \mu_{i+j+r-2}(q)\bigr)_{1\leq i,j\leq2n},
\]
where $r$ is a fixed integer.
We mainly investigate the case where $\mu_{n}$ is
the moments of the little $q$-Jacobi polynomials in this paper.
%% We first recall some terminology in $q$-series.
Throughout this paper we use the standard  notation for $q$-series (see \cite{AAR,GR}):
\begin{equation*}
    (a;q)_{\infty}=\prod_{k=0}^{\infty}(1-aq^{k}),\qquad
    (a;q)_{n}=\frac{(a;q)_{\infty}}{(aq^{n};q)_{\infty}}
\end{equation*}
for any integer $n$.
Usually  $(a;q)_{n}$ is called the  \defterm{$q$-shifted factorial},
and we frequently use the compact notation:
\begin{align*}
    &(a_{1},a_{2},\dots,a_{r};q)_{\infty}=(a_{1};q)_{\infty}(a_{2};q)_{\infty}\cdots(a_{r};q)_{\infty},\\
    &(a_{1},a_{2},\dots,a_{r};q)_{n}=(a_{1};q)_{n}(a_{2};q)_{n}\cdots(a_{r};q)_{n}.
\end{align*}
The \defterm{${}_{r+1}\phi_{r} $ basic hypergeometric series} is defined by
\begin{align*}
{}_{r+1}\phi_{r}\left[\,
{{a_{1},a_{2},\dots,a_{r+1}}\atop{b_{1},\dots,b_{r}}};q,z
\,\right]
=\sum_{n=0}^{\infty}\frac{(a_{1},a_{2},\dots,a_{r+1};q)_{n}}{(q,b_{1},\dots,b_{r};q)_{n}}z^{n}.
\end{align*}
The \defterm{little $q$-Jacobi polynomials} \cite{GR,KLS} are defined by
\begin{equation}
p_{n}(x;a,b;q)=\frac{(aq;q)_{n}}{(abq^{n+1};q)_{n}}(-1)^{n}q^{\binom{n}2}
{}_{2}\phi_{1}\left[
{{q^{-n}, abq^{n+1}}
\atop{aq}}
\,;\, q, xq
\right],
\label{eq:little-q-Jacobi}
\end{equation}
which are orthogonal with respect to the \defterm{inner product} defined by
\begin{equation}
\langle f,g\rangle=\frac{(aq;q)_{\infty}}{(abq^2;q)_{\infty}}
\sum_{k=0}^{\infty}
\frac{(bq;q)_{k}}{(q;q)_{k}}(aq)^{k}f\left(q^{k}\right)g\left(q^{k}\right).
\label{eq:inner-product}
\end{equation}
The moments of the little $q$-Jacobi polynomials are defined by
\[
\mu_{n}=\langle x^{n},1\rangle=\frac{(aq;q)_{n}}{(abq^2;q)_{n}}\quad(n=0,1,2,\dots).
\]
The main results on 
$
\Pf\Bigl((q^{i-1}-q^{j-1})\mu_{i+j-2}\Bigr)_{1\leq i,j\leq 2n}
$
 are stated in Section~\ref{sec:q-Catalan-Pfaffians}.
%% If we put $a=q^{\alpha}$, $b=q^{\beta}$ and let $q\rightarrow1$,
%% then we obtain
%% \[
%%   \mu_{n}\rightarrow\frac{(\alpha+1)_{n}}{(\alpha+\beta+2)_{n}},
%% \]
%% hence various Catalan-like numbers are obtained as speciaizations of the moments as follows:
%% \[
%% \frac{(\frac12)_{n}}{(2)_{n}}=\frac{C_{n}}{2^{2n}},\quad
%% \frac{(\frac12)_{n}}{(1)_{n}}=\frac{1}{2^{2n}}\binom{2n}{n},\quad
%% \frac{(\frac32)_{n}}{(2)_{n}}=\frac{1}{2^{2n}}\binom{2n+1}{n}.
%% \]
%% This is the reason we take the moments as a $q$-analogue of the Catalan numbers $C_{n}$.

To prove the Pfaffian identities
we employ an $LU$-type decomposition of a skew-symmetric matrix,
which we call a Pfaffian decomposition.
In Section~\ref{sec:decomposition} we state this decomposition and give a proof
by using a Pfaffian analogue of the Desnanot-Jacobi adjoint matrix theorem \cite[Theorem~3.12]{B}.
%% In the next section we consider a decomposition for skew-symmetric matrices.
%% The $LU$-decomposition is a well-known method to write any square matrix $A$ in the form
%% $PA=LDU$ where $P$ is a permutation matrix, $D$ a diagonal matrix,
%% $L$ (resp. $U$) a lower (resp. upper) unitriangular matrix.
%% Usually we can take the identity matrix for $P$.
%% We give a proof to make this paper comprehensive and show that 
%% each entry of $D$, $L$ and $U$ is written by the determinants of certain type of sub-matrices $A$
%% using the Desnanot-Jacobi adjoint matrix theorem \cite[Theorem~3.12]{B}.
%% This means the following.
%
In Section~\ref{sec:proof}
we give a proof of our main results stated in Section~\ref{sec:q-Catalan-Pfaffians}.
We prove the Pfaffian decomposition in Theorem~\ref{conj:decomposition}
by reducing the single sum obtained as the matrix multiplication to the $q$-Dougall formula \eqref{eq:q-Dougall}
for a terminating very-well-poised ${}_{6}\phi_{5}$ series
(see \cite{AAR,GR}).
As a byproduct of the proof
we obtain the Pfaffian decomposition of another skew-symmetric matrix
stated in Theorem~\ref{th:byproduct}.
%% Assume we are given a square matrix $A$ of size $n$ which we don't know
%% how to evaluate the determinant.
%% Suppose we can take the identity for the permutation matrix $P$,
%% that is to say, $\det A^{[k]}_{[k]}\neq0$ for $1\leq k\leq n$.
%% If we could guess a formula for $\det A^{[i]}_{[i-1],j}$ 
%% and $\det A^{[j-1],i}_{[j]}$ for $1\leq i,j\leq n$,
%% then we would obtain an identity to prove which usually contain only a single sum.
\smallskip

As an application of our main results in Section~\ref{sec:q-Catalan-Pfaffians} and Section~\ref{sec:proof},
we obtain a formula for weighted enumeration of shifted reverse plane partitions
in Section~\ref{sec:application}.
We consider a special family of shifted reverse plane partitions
and give weights that resembles to the weight in the inner product \eqref{eq:inner-product}
to profiles of shifted reverse plane partitions in the family
(see \eqref{eq:profile-weight} and \eqref{eq:profile-weight2}).
Then Corollary~\ref{conj:01} (resp. Corollary~\ref{cor:byproduct})
gives the weighted enumeration of the family of shifted reverse plane partitions
whose number of rows are even (resp. odd).
\smallskip 

In Section~\ref{sec:open-problems}
we state several conjectures for this type of Pfaffians.
One may ask what is the relation between our Pfaffians and the classical theory
of orthogonal polynomials.
At this point we have no answer to the question why the Pfaffians factors into nice linear factors
from the view point of the classical theory.
\smallskip

Finally, in Appendix we state our second proof of the main results
in Section~\ref{sec:proof}
using  Zeilberger's creative telescoping.
We see that the certificates are simple,
and we can prove the formulas by hand.

%% file: tagawa02.tex
%%%% ---------------< ooo >--------------- %%%%
%%%%
%%%% Section 2
%%%%
%%%% ---------------< ooo >--------------- %%%%

\section{Pfaffian decomposition}
\label{sec:decomposition}

First we recall the reader a well-known decomposition of a matrix.
Let $A=(a^{i}_{j})_{i,j\geq1}$ be a matrix (of finite or infinite row/column length).
If $I=\left\{i_{1},\dots,i_{r}\right\}$ (resp. $J=\left\{j_{1},\dots,j_{r}\right\}$) 
are a set of row (resp. column) indices,
then we write $A^{I}_{J}=A^{i_{1},\dots,i_{r}}_{j_{1},\dots,j_{r}}$
for the $r\times r$ submatrix obtained from $A$ by choosing 
the rows indexed by $I$ and columns indexed by $J$.
Let $a^{I}_{J}=a^{i_{1},\dots,i_{r}}_{j_{1},\dots,j_{r}}$ denote $\det A^{I}_{J}$ if $|I|=|J|>0$,
and $1$ if $I=J=\emptyset$.
%% For positive integer $n$ we let $[n]=\{1,\ldots,n\}$.
The following identity is known as the Desnanot-Jacobi adjoint matrix theorem \cite[Theorem~3.12]{B}
%
%
%% ------------------------- %%
%%  Desnanot-Jacobi
%% ------------------------- %%
%
%
\begin{eqnarray}
&&
\det A^{[n-2]}_{[n-2]}
\det A^{[n]}_{[n]}
\nonumber\\
&&\qquad=\det A^{[n-2],n-1}_{[n-2],n-1}
\det A^{[n-2],n}_{[n-2],n}
-\det A^{[n-2],n-1}_{[n-2],n}
\det A^{[n-2],n}_{[n-2],n-1}.
\label{eq:Desnanot-Jacobi}
\end{eqnarray}
The following proposition is usually called the $LU$-decomposition of a matrix.
Usually $LU$-decomposition means $L$ is lower unitriangular and $U$ is upper triangular.
But here we adopt the style of $LDU$-decomposition where both of $L$ and $U$
are (lower or upper) unitriangular and $D$ is diagonal
because of Theorem~\ref{th:Pfaffian-decomp}.
If $A$ is an $n\times n$ matrix of rank $n$,
then, by elementary linear algebra,
we can deduce that there is an $n\times n$ permutation matrix $P$ such that
$\det \left(PA\right)^{[i]}_{[i]}\neq0$ for any $i=1,\dots,n$.
Although this permutation matrix $P$ is not unique,
the triple $(L,D,U)$ in the following theorem is unique for chosen $A$ and $P$.
But this is the most general case, and in many applications we can choose $P$ to be the identity
matrix $I_n$.
We give a proof here to make this paper more comprehensive
and as a warm-up for the succeeding proof of Theorem~\ref{th:Pfaffian-decomp}.
%
%
%% ------------------------- %%
%% Theorem
%% ------------------------- %%
%
%
\begin{prop}
Let $n$ be a positive integer, 
and $A=(a^{i}_{j})_{1\leq i,j\leq n}$ be an $n\times n$ matrix of rank $n$.
If we choose an $n\times n$ permutation matrix $P$ 
such that $\det \left(PA\right)^{[i]}_{[i]}\neq0$ for $1\leq i\leq n$,
then $PA$ is uniquely written as
\begin{equation}
PA=L\,D\,U,
\label{eq:LU-decomposition}
\end{equation}
where $D=(d_{i}\delta^{i}_{j})_{1\leq i,j\leq n}$ is a diagonal matrix,
$L=(l^{i}_{j})_{1\leq i,j\leq n}$ is a lower unitriangular matrix
and $U=(u^{i}_{j})_{1\leq i,j\leq n}$ is an upper unitriangular matrix.
In fact
\begin{eqnarray*}
d_{i}=\frac{\det \left(PA\right)_{[i]}^{[i]}}{\det \left(PA\right)_{[i-1]}^{[i-1]}},
\qquad
l^{i}_{j}=\frac{\det \left(PA\right)_{[j]}^{[j-1],i}}{\det \left(PA\right)_{[j]}^{[j]}},
\qquad
u^{i}_{j}=\frac{\det \left(PA\right)^{[i]}_{[i-1],j}}{\det \left(PA\right)_{[i]}^{[i]}}.
\end{eqnarray*}
Here the Kronecker delta $\delta^{i}_{j}$
takes the value 1 if $i=j$, and $0$ otherwise.
%% =\begin{cases}1,&\text{ if $i=j$}\\0,&\text{ if $i\neq j$,}\end{cases}$ is the Kronecker delta.
\end{prop}
%
%
%% ------------------------- %%
%% Proof of Theorem
%% ------------------------- %%
%
%
\begin{demo}{Proof}
We may assume $P=I_{n}$ without loss of generality,
and we have to show that
$a^{i}_{j}$ is uniquely written as
$
a^{i}_{j}
=\sum_{k=1}^{\min(i,j)}
l^{i}_{k}d_{k}u^{k}_{j}
$
with $l^{i}_{i}=u^{i}_{i}=1$ for $1\leq i\leq n$.
This is trivial if $n=1$.
Assume $n\geq2$ and this is true for all $1\leq i,j\leq n-1$.
If $1\leq i< n$,
then $u^{i}_{n}$ must satisfy
\[
a^{i}_{n}=\sum_{k=1}^{i-1}
\frac{a^{[k-1],i}_{[k]}\,a^{[k]}_{[k-1],n}}{a^{[k-1]}_{[k-1]}\,a^{[k]}_{[k]}}
+d_{i}u^{i}_{n}.
\]
But this can be obtained from 
\[
a^{i}_{i}=\sum_{k=1}^{i-1}
\frac{a^{[k-1],i}_{[k]}\,a^{[k]}_{[k-1],i}}{a^{[k-1]}_{[k-1]}\,a^{[k]}_{[k]}}
+\frac{a^{[i]}_{[i]}}{a^{[i-1]}_{[i-1]}}
\]
by replacing $i$th column by $n$th column of $A$,
and we obtain
$d_{i}u^{i}_{n}={a^{[i]}_{[i-1],n}}/{a^{[i-1]}_{[i-1]}}$.
Hence we derive $u^{i}_{n}={a^{[i]}_{[i-1],n}}/{a^{[i]}_{[i]}}$,
and vice versa.
Similarly,
when $1\leq j\leq n-1$,
we can show that $l^{n}_{j}$ is determined uniquely and given by the above formula.
Thus it is enough to prove the formula for $i=j=n$,
which implies
\[
a^{n}_{n}=\sum_{k=1}^{n-1}
\frac{a^{[k-1],n}_{[k]}\,a^{[k]}_{[k-1],n}}{a^{[k-1]}_{[k-1]}\,a^{[k]}_{[k]}}
+d_{n}.
\]
By induction hypothesis
\[
a^{n-1}_{n-1}
=\sum_{k=1}^{n-2}
\frac{a^{[k-1],n-1}_{[k]}\,a^{[k]}_{[k-1],n-1}}{a^{[k-1]}_{[k-1]}\,a^{[k]}_{[k]}}
+\frac{a^{[n-1]}_{[n-1]}}{a^{[n-2]}_{[n-2]}},
\]
which implies
\[
a^{n}_{n}
-\sum_{k=1}^{n-2}
\frac{a^{[k-1],n}_{[k]}\,a^{[k]}_{[k-1],n}}{a^{[k-1]}_{[k-1]}\,a^{[k]}_{[k]}}
=\frac{a^{[n-2],n}_{[n-2],n}}{a^{[n-2]}_{[n-2]}}.
\]
Hence
\[
d_{n}
=a^{n}_{n}
-\sum_{k=1}^{n-1}
\frac{a^{[k-1],n}_{[k]}\,a^{[k]}_{[k-1],n}}{a^{[k-1]}_{[k-1]}\,a^{[k]}_{[k]}}
=\frac{a^{[n-2],n}_{[n-2],n}}{a^{[n-2]}_{[n-2]}}
-\frac{a^{[n-2],n}_{[n-1]}\,a^{[n-1]}_{[n-2],n}}{a^{[n-2]}_{[n-2]}\,a^{[n-1]}_{[n-1]}},
\]
which equals ${a^{[n]}_{[n]}}/{a^{[n-1]}_{[n-1]}}$
by \eqref{eq:Desnanot-Jacobi}.
Conversely, if we take $d_{n}={a^{[n]}_{[n]}}/{a^{[n-1]}_{[n-1]}}$,
then it clearly satisfies the above equation and gives the $LU$-decomposition of $A$.
\end{demo}
%
%
%
%% -------------------------------------------------- %%
%% Pfaffian decomposition
%% -------------------------------------------------- %%
The fact that each entry of $L$, $D$ and $U$ is expressed with certain type of minors of $A$
 appears in \cite{N} related to the Painlev\'e equations.
Although we can use this decomposition even in the case where $A$ is a skew-symmetric matrix,
it seems more consistent to consider the decomposition in the following theorem
in which each entry is expressed with subpfaffians.
This type of decomposition also seems important with relation to the integrable systems
(see \cite{AM}).
Let us start with some definitions.

\smallskip

We define $2\times2$ skew-symmetric matrix $J_{2}$ by
\[
J_{2}=\begin{pmatrix}
0&1\\
-1&0
\end{pmatrix},
\]
and let $J_{2n}=J_{2}\oplus\dots\oplus J_{2}$ denote the $2n\times2n$ matrix
whose main diagonal $2\times2$ blocks are all $J_{2}$ and the other blocks
are $2\times2$ zero matrices $O_{2}$.
Note that ${}^t\!J_{2n} J_{2n}=I_{2n}$.

For a skew-symmetric matrix $A$,
we usually take $I=J$
so, hereafter,
we write $A_{I}=A_{i_{1},\dots,i_{r}}$ for $A^{I}_{I}$.
Further let $a_{I}=a_{i_{1},\dots,i_{r}}$
denote $\Pf A_{I}$ if $I\neq\emptyset$, $1$ if $I=\emptyset$ when there is no fear of confusion.
%
%
%% ------------------------- %%
%%  Desnanot-Jacobi (Pfaffian)
%% ------------------------- %%
%
%
%% In this paper we call the following identity
Then the Pfaffian analogue of the Desnanot-Jacobi adjoint-matrix theorem \eqref{eq:Desnanot-Jacobi}
 reads as follows
(see  \cite{IW,Kn}):
\begin{eqnarray}
&&
a_{[n-4]}
a_{[n]}
=a_{[n-4],n-3,n-2}
a_{[n-4],n-1,n}
\label{eq:Pf-Desnanot-Jacobi}
\\&&\qquad
-a_{[n-4],n-3,n-1}
a_{[n-4],n-2,n}
+a_{[n-4],n-3,n}
a_{[n-4],n-2,n-1}.
\nonumber
\end{eqnarray}
The following theorem gives 
so-called \defterm{Pfaffian decomposition} of a skew-symmetric matrix $A$.
%%which gives a decomposition of any skew-symmetric matrix.
%
%
%% ------------------------- %%
%% Theorem
%% ------------------------- %%
%
%
\begin{theorem}
\label{th:Pfaffian-decomp}
Let $n$ be a positive integer, and $A=(a^{i}_{j})_{1\leq i,j\leq2n}$ be a skew-symmetric matrix of size $2n$.
If $a_{[2i]}\neq0$ for $1\leq i\leq n$,
then $A$ is uniquely written as
\begin{equation}
A={}^{t}V\,T\,V.
\label{eq:Pfaffian-decomposition}
\end{equation}
Here $T$ and $V$ are composed of $2\times2$ blocks
\begin{align*}
T=\left(
\begin{matrix}
T_{1}&O_{2}&\hdots&O_{2}\\
O_{2}&T_{2}&\hdots&O_{2}\\
\vdots&\vdots&\ddots&\vdots\\
O_{2}&O_{2}&\hdots&T_{n}
\end{matrix}
\right),
\qquad
V=\left(
\begin{matrix}
J_{2}&V^{1}_{2}&\hdots&V^{1}_{n}\\
O_{2}&J_{2}&\hdots&V^{2}_{n}\\
\vdots&\vdots&\ddots&\vdots\\
O_{2}&O_{2}&\hdots&J_{2}
\end{matrix}
\right),
\end{align*}
of the form
$
T_{i}=\begin{pmatrix}
0&t_{i}\\
-t_{i}&0
\end{pmatrix}
$
for $1\leq i\leq n$,
and
$
V^{i}_{j}=\begin{pmatrix}
v^{2i-1}_{2j-1}(i)&v^{2i-1}_{2j}(i)\\
v^{2i}_{2j-1}(i)&v^{2i}_{2j}(i)
\end{pmatrix}
$
for $1\leq i<j\leq n$,
where each $t_{i}$ and $v^{k}_{l}(i)$ is defined by
\begin{eqnarray}
t_{i}=\frac{a_{[2i]}}{a_{[2i-2]}},
\qquad\qquad
v^{k}_{l}(i)=\frac{a_{[2i-2],k,l}}{a_{[2i]}}
\label{eq:T_and_U}
\end{eqnarray}
for $1\leq i\leq n$ and $1\leq k,l\leq 2n$.
\end{theorem}
Before we proceed to the proof of the theorem,
we illustrate the decomposition by an example.
If we take a $4\times4$ skew-symmetric matrix $A=(a_{ij})_{1\leq i,j\leq4}$,
then the above decomposition is given by
\begin{align*}
&T=\begin{pmatrix}
0&a_{12}&0&0\\
-a_{12}&0&0&0\\
0&0&0&\frac{a_{1234}}{a_{12}}\\
0&0&-\frac{a_{1234}}{a_{12}}&0\\
\end{pmatrix},\qquad
V=\begin{pmatrix}
0&1&\frac{a_{13}}{a_{12}}&\frac{a_{14}}{a_{12}}\\
-1&0&\frac{a_{23}}{a_{12}}&\frac{a_{24}}{a_{12}}\\
0&0&0&1\\
0&0&-1&0\\
\end{pmatrix},
\end{align*}
where $a_{ij}=\Pf\begin{pmatrix}0&a^{i}_{j}\\-a^{i}_{j}&0\end{pmatrix}=a^{i}_{j}$
and $a_{1234}=\Pf A$.
%
%
%% ------------------------- %%
%% Proof of Theorem
%% ------------------------- %%
%
%
\begin{demo}{Proof of Theorem~\ref{th:Pfaffian-decomp}}
First we write the matrix $A$ by $2\times2$ blocks as
\[
A=\begin{pmatrix}
A^{1}_{1}&A^{1}_{2}&\hdots&A^{1}_{n}\\
A^{2}_{1}&A^{2}_{2}&\hdots&A^{2}_{n}\\
\vdots&\vdots&\ddots&\vdots\\
A^{n}_{1}&A^{n}_{2}&\hdots&A^{n}_{n}
\end{pmatrix},
\]
where $A^{i}_{j}$ is the $2\times2$ block matrix
$
A^{i}_{j}=\begin{pmatrix}
a^{2i-1}_{2j-1}&a^{2i-1}_{2j}\\
a^{2i}_{2j-1}&a^{2i}_{2j}
\end{pmatrix}
$
for $1\leq i,j\leq n$.
Hence the decomposition \eqref{eq:Pfaffian-decomposition} is equivalent to
\begin{equation}
A^{i}_{j}=\sum_{k=1}^{\min(i,j)}{}^{t}V^{k}_{i}\,T_{k}\,V^{k}_{j}
\label{eq:claim}
\end{equation}
with $V^{i}_{i}=J_{2}$.
We proceed by induction on $n$.
If $n=1$,
then \eqref{eq:claim} implies $T=T_{1}=A^{1}_{1}=A$
and $V=V^{1}_{1}=J_{2}$
so that the existence and uniqueness are trivial.
Assume $n\geq2$,
 and our claim holds for $n-1$.
That is, the equations \eqref{eq:claim} 
for $1\leq i,j<n$ uniquely determines
all $T_{i}$ and $V^{i}_{j}$ for $1\leq i,j< n$
and each entry is given by \eqref{eq:T_and_U}.
This implies that
\begin{align*}
\sum_{k=1}^{i-1}\left(
\frac{a_{[2k-1],2i-1}a_{[2k-2],2k,2i}}{a_{[2k-2]}a_{[2k]}}
-\frac{a_{[2k-2],2k,2i-1}a_{[2k-1],2i}}{a_{[2k-2]}a_{[2k]}}
\right)
+\frac{a_{[2i]}}{a_{[2i-2]}}
=a^{2i-1}_{2i}
\end{align*}
holds for $1\leq i< n$ from \eqref{eq:claim}.
Replacing $(2i-1)$st row/column by $r$th row/column
and $2i$th row/column by $s$th row/column in this identity,
we see
\begin{align}
\sum_{k=1}^{i-1}\left(
\frac{a_{[2k-1],r}a_{[2k-2],2k,s}}{a_{[2k-2]}a_{[2k]}}
-\frac{a_{[2k-2],2k,r}a_{[2k-1],s}}{a_{[2k-2]}a_{[2k]}}
\right)
+\frac{a_{[2i-2],r,s}}{a_{[2i-2]}}
=a^{r}_{s}
\label{eq:plucker}
\end{align}
holds for any $r$ and $s$.
From computation of each entry of the equation \eqref{eq:claim},
we see that $v^{r}_{s}$ ($1\leq i<n$, $r=2i-1,2i$, $s=2n-1,2n$) must satisfy
\[
\sum_{k=1}^{i-1}\left(
\frac{a_{[2k-1],r}a_{[2k-2],2k,s}}{a_{[2k-2]}a_{[2k]}}
-\frac{a_{[2k-2],2k,r}a_{[2k-1],s}}{a_{[2k-2]}a_{[2k]}}
\right)
+\frac{a_{[2i]}}{a_{[2i-2]}}v^{r}_{s}
=a^{r}_{s}.
\]
Comparing this equation with \eqref{eq:plucker},
we see that $v^{r}_{s}$ in \eqref{eq:T_and_U} is the unique solution of this equation.
Similarly,
from computation of each entry of the equation \eqref{eq:claim},
$t_{n}$ must satisfy
\[
\sum_{k=1}^{n-1}\left(
\frac{a_{[2k-1],2n-1}a_{[2k-2],2k,2n}}{a_{[2k-2]}a_{[2k]}}
-\frac{a_{[2k-2],2k,2n-1}a_{[2k-1],2n}}{a_{[2k-2]}a_{[2k]}}
\right)
+t_{n}
=a^{2n-1}_{2n}.
\]
Substituting $i=n-1$, $r=2n-1$ and $s=2n$ into \eqref{eq:plucker},
we obtain
\[
\sum_{k=1}^{n-2}\left(
\frac{a_{[2k-1],2n-1}a_{[2k-2],2k,2n}}{a_{[2k-2]}a_{[2k]}}
-\frac{a_{[2k-2],2k,2n-1}a_{[2k-1],2n}}{a_{[2k-2]}a_{[2k]}}
\right)
+\frac{a_{[2n-4],2n-1,2n}}{a_{[2n-4]}}
=a^{2n-1}_{2n}.
\]
Hence we have
\[
t_{n}=\frac{a_{[2n-4],2n-1,2n}}{a_{[2n-4]}}
-\frac{a_{[2n-3],2n-1}a_{[2n-4],2n-2,2n}}{a_{[2n-4]}a_{[2n-2]}}
+\frac{a_{[2n-4],2n-2,2n-1}a_{[2n-3],2n}}{a_{[2n-4]}a_{[2n-2]}}.
\]
Thus, by \eqref{eq:Pf-Desnanot-Jacobi},
we conclude that $t_n$ in \eqref{eq:T_and_U} is the unique solution of this equation,
and this proves the theorem in the case of $n$.
\end{demo}
This theorem shows that, if we obtain a guess for each entry of $T$ and $V$,
then, by uniqueness of the decomposition, it is enough to prove the matrix multiplication,
 which is equivalent to the single sum
\begin{align*}
\sum_{k\geq1}\left\{
v_{i}^{2k-1}(k)t_{k}v_{j}^{2k}(k)
-v_{i}^{2k}(k)t_{k}v_{j}^{2k-1}(k)
\right\}
=a^{i}_{j}.
\end{align*}
From \eqref{eq:T_and_U}
it is enough to guess a formula for the subpfaffians $a_{[2i-1],j}$
and $a_{[2i-2],2i,j}$ for any row/column indices $i$ and $j$.

\smallskip
%
%
%% ------------------------- %%
%% LU-decomposition
%% ------------------------- %%
%
%
%
By the uniqueness of $LU$-decomposition \eqref{eq:LU-decomposition}
and
Pfaffian decomposition \eqref{eq:Pfaffian-decomposition},
the $LU$-decomposition and the Pfaffian decomposition are, in a sense,
equivalent.
We can get the $LU$-decomposition from the Pfaffian decomposition, and vise versa.
%Assume $\Pf A_{[2i]}\neq0$ for any $i\geq1$.
If we put $P=(p^{i}_{j})_{i,j\geq1}$,
where
\[
p^{i}_{j}
=\begin{cases}
\delta^{i+1}_{j}
&\text{ if $i$ is odd,}\\
\delta^{i-1}_{j}
&\text{ if $i$ is even,}
\end{cases}
\]
which is the permutation matrix corresponding to $(12)(34)\dots$,
then it is easy to see that $\det\left(PA\right)^{[i]}_{[i]}\neq0$
for $i\geq1$.
If we put 
$\displaystyle
J=\bigoplus_{i\geq1} J_{2}
$,
then
$U={}^{t}\!J\,V$ is upper unitriangular,
$D=P\,T=P\,{}^{t}\!J\,T\,J$ is diagonal,
$L=P\,{}^tV\,J\,P$ is lower unitriangular,
hence
\begin{equation}
P\,A=L\,D\,U
\end{equation}
gives the $LU$-decomposition.
Each entry of the matrices
$U=\left(u^{i}_{j}\right)_{i,j\geq 1}$,
$L=\left(l^{i}_{j}\right)_{i,j\geq 1}$
and $D=\left(d_{i}\,\delta^{i}_{j}\right)_{i,j\geq 1}$
is given by
\begin{align*}
u^{i}_{j}
&=\begin{cases}
-v^{i+1}_{j}
&\text{ if $i$ is odd,}\\
v^{i-1}_{j}
&\text{ if $i$ is even,}
\end{cases}
\qquad\qquad
d_{i}=\begin{cases}
-t_{(i+1)/2}
&\text{ if $i$ is odd,}\\
t_{i/2}
&\text{ if $i$ is even,}
\end{cases}
\\
l^{i}_{j}
&=\begin{cases}
v^{j}_{i+1}
&\text{ if $i$ is odd and $j$ is odd,}\\
v^{j}_{i-1}
&\text{ if $i$ is even and $j$ is odd,}\\
-v^{j}_{i+1}
&\text{ if $i$ is odd and $j$ is even,}\\
-v^{j}_{i-1}
&\text{ if $i$ is even and $j$ is even.}
\end{cases}
\end{align*}
\smallskip
For later use
we cite the minor summation formula of Pfaffians here:
%%%%%%%%%%%%%%%%%%%%%%%%%%%%%%%%%%%%%%%%%%%%%%%%%%%%
%
% Minor Summation Formula
%
%%%%%%%%%%%%%%%%%%%%%%%%%%%%%%%%%%%%%%%%%%%%%%%%%%%%
\begin{theorem}
\label{th:msf}
(\cite{IW1,IW})
Let $n\le N$ be positive integers and assume $n$ is even.
Let $T=(t^{i}_{j})_{1\le i\le n, 1\le j\le N}$ be an $n\times N$ rectangular matrix,
and let $B=(b^{i}_{j})_{1\le i,j\le N}$ be a skew symmetric matrix of size $N$.
Then we have
\begin{align}
\label{eq_msf}
\sum_{{I\subseteq[N]}\atop{\sharp I=n}}
\Pf(B_I) \det(T^{[n]}_I)&
=\Pf(Q),
%% &=\Pf(A)^{-M+m+1}\Pf\begin{pmatrix}
%%    O_n   &   T J_N \\
%% -J_N {}^tT & \widehat{A} \\
%% \end{pmatrix},\nonumber\\
\end{align}
where the skew symmetric matrix $Q=(Q^{i}_{j})=TB\,{}^t\kern-1pt T$ of size $n$ 
whose entries are given by
\begin{equation}
\label{pf_msf}
Q^{i}_{j}=\sum_{1\le k<l\le N} b^{k}_{l} \det(T^{ij}_{kl}),
\qquad(1\le i,j\le n).
\end{equation}
%Here we write $T^{ij}_{kl}$ for $T^{\{ij\}}_{\{kl\}}=\begin{vmatrix}t_{ik}&t_{il}\\t_{jk}&t_{jl}\end{vmatrix}$.
\end{theorem}
When $n$ is odd,
we can immediately derive a similar formula from the case when $n$ is even.
%
%
%
%
%
%% ------------------------- %%
%% Proposition for Subpfaffian
%% ------------------------- %%
%
%
\begin{prop}
\label{prop:subpfaffian}
Let $\{\alpha_{k}\}_{k\geq1}$ be any sequence,
and let $n$ be a positive integer.
Set $B=(b^{i}_{j})_{i,j\geq1}$ to be the skew-symmetric matrix defined by
\begin{equation}
b^{i}_{j}=\begin{cases}
\alpha_{i}&\text{ if $j=i+1$ for $i\geq1$,}\\
-\alpha_{j}&\text{ if $i=j+1$ for $j\geq1$,}\\
0&\text{ otherwise.}
\end{cases}
\label{eq:alpha_k}
\end{equation}
If $I=(i_{1},\dots,i_{2n})$ is an index set
such that $1\leq i_{1}<\dots<i_{2n}$,
then
\begin{equation}
\Pf\left(B_{I}\right)
=\begin{cases}
\prod_{k=1}^{n} \alpha_{i_{2k-1}}
&\text{ if $i_{2k}=i_{2k-1}+1$ for $k=1,\dots,n$,}\\
0&\text{ otherwise.}
\end{cases}
\end{equation}
\end{prop}
%
%
%

%% file: tagawa03.tex
%%%% ---------------< ooo >--------------- %%%%
%%%%
%%%% Section 3
%%%%
%%%% ---------------< ooo >--------------- %%%%
%
%
%
\section{A Pfaffian analogue of $q$-Catalan Hankel determinants}
\label{sec:q-Catalan-Pfaffians}
%
%
%
%We use the notation
%\begin{equation*}
%(a;q)_{n}=\begin{cases}
%\prod_{i=0}^{n-1}(1-aq^{i})
%&\text{ if $n\geq0$,}\\
%1/\prod_{i=0}^{-n-1}(1-aq^{i+n})
%&\text{ if $n<0$.}\\
%\end{cases}
%\end{equation*}
%
%
%
%
%
Let us write
\begin{equation}
a^{i}_{j}=(q^{i-1}-q^{j-1})\frac{(aq;q)_{i+j+r-2}}{(abq^2;q)_{i+j+r-2}}
\label{eq:a^{i}_{j}}
\end{equation}
for $i,j\geq1$,
and let $A$ denote the skew-symmetric matrix
\[
A=\left(a^{i}_{j}\right)_{i,j\geq 1}
\]
of infinite degree.
Then the following theorem gives the Pfaffian decomposition of $A$.
%
%
%
%
%
%% ------------------------- %%
%% Conjecture
%% ------------------------- %%
%
%
\begin{theorem}
\label{conj:decomposition}
Let $A$ be as above,
and let
\begin{align*}
t_{i}
&=a^{i-1}q^{(i-1)(i+r)}
\frac{
(q;q)_{i}(aq;q)_{i+r}(bq;q)_{i-1}
}
{
(abq^2;q)_{2i+r-1}
(abq^{i+r};q)_{i-1}
},\\
v^{i}_{j}
&=\begin{cases}
o^{i}_{j}&\text{ if $i$ is odd,}\\
e^{i}_{j}&\text{ if $i$ is even,}
\end{cases}
\end{align*}
where
\begin{align*}
o^{i}_{j}
&=\frac{
(q^{j-i};q)_{i}
}{
(q;q)_{i}
}
\cdot\frac{
(aq^{i+r+1};q)_{j-i-1}
}
{
(abq^{2i+r+1};q)_{j-i-1}
},
\\
e^{i}_{j}
&=
q\frac{
(q^{j-i};q)_{1}
(q^{j-i+2};q)_{i-2}
}{
(q;q)_{i-1}
}
\cdot\frac{
(aq^{i+r};q)_{j-i}f(i,j,r)
}
{
(abq^{2i+r-3};q)_{1}
(abq^{2i+r-1};q)_{j-i+1}
},
\end{align*}
with
\begin{align}
&f(i,j,r) 
%&= [i-1]_q 
%-  aq^{i+r-1} \Bigl\{[j]_q - q^{i-2}[j-i+1]_q
%\nonumber\\&\qquad\qquad
%+ bq^{i-2} ([j]_q - q[j - i+1]_q) \Bigr\}
%+ a^2bq^{2i+j+2r-3}[i-1]_q
%\nonumber\\
=(1-q^{i-1})(1-aq^{i+r-1})(1-abq^{i+j+r-2})/(1-q)
\nonumber\\&\qquad\qquad\qquad\qquad\qquad
+aq^{2i+r-3}(1-b)(1-q^{j-i+1}).
\label{eq:f(i,j)}
\end{align}
If we put
$\displaystyle
T=\bigoplus_{i\geq1}\begin{pmatrix} 0&t_{2i-1}\\-t_{2i-1}&0\end{pmatrix}
$
and
$
V=\left(v^{i}_{j}\right)_{i,j\geq1}
$
then
\begin{equation}
A={}^{t}V\, T\, V
\label{eq:decomposition}
\end{equation}
gives the Pfaffian decomposition of $A$.
\end{theorem}
An immediate consequence of the theorem is the following corollary.
%
%
%% ------------------------- %%
%% Catalan Pfaffian
%% ------------------------- %%
%
%
\begin{corollary}
\label{conj:01}
Let $n\geq1$ and $r$ be integers.
Then we have
\begin{align}
&\Pf\left((q^{i-1}-q^{j-1})\frac{(aq;q)_{i+j+r-2}}{(abq^2;q)_{i+j+r-2}}\right)_{1\leq i,j\leq 2n}
\nonumber\\&=
a^{n(n-1)}q^{n(n-1)(4n+1)/3+n(n-1)r}\prod_{k=1}^{n-1}(bq;q)_{2k}
\prod_{k=1}^n\frac{(q;q)_{2k-1}(aq;q)_{2k+r-1}}
{(abq^2;q)_{2(k+n)+r-3}}.
\label{eq:pf-special}
\end{align}
%\end{corollary}
%
%
%
%
In fact, we obtain a more general formula from Theorem~\ref{conj:decomposition}.
%
%
%% ------------------------- %%
%% Conjecture
%% ------------------------- %%
%
%
%\begin{corollary}
%Let $n\geq1$ and $m\geq2n$ be integers.
If $A$ is as above and $m$ is a positive integer,
then the following identities hold:
\begin{align}
&\Pf\left(A_{[2n-1],m}\right)
=a^{n(n-1)}q^{n(n-1)(4n+1)/3+n(n-1)r}
\nonumber\\&\times
\frac{(q^{m-2n+1};q)_{2n-1}(aq;q)_{m+r-1}}
{(abq^2;q)_{m+2n+r-3}}
\prod_{k=1}^{n-1}\frac{
(bq;q)_{2k}
(q;q)_{2k-1}
(aq;q)_{2k+r-1}
}
{(abq^2;q)_{2(k+n)+r-3}},
\label{eq:pf-general1}
\\
&\Pf\left(A_{[2n-2],2n,m}\right)
=a^{n(n-1)}q^{n(n-1)(4n+1)/3+n(n-1)r+1}
f(2n,m,r)
\nonumber\\&\times
\frac{(q^{m-2n};q)_{1}(q^{m-2n+2};q)_{2n-2}(aq;q)_{m+r-1}}
{(abq^{4n+r-3};q)_{1}(abq^2;q)_{m+2n+r-2}}
\prod_{k=1}^{n-1}\frac{
(bq;q)_{2k}
(q;q)_{2k-1}
(aq;q)_{2k+r-1}
}
{(abq^2;q)_{2(k+n)+r-3}},
\label{eq:pf-general2}
\end{align}
where $f(i,j,r)$ is defined by \eqref{eq:f(i,j)}.
\end{corollary}

\smallskip
Next we consider a specialization of Corollary~\ref{conj:01}.
If we put $a=q^{\alpha}$ and $b=q^{\beta}$
and let $q\rightarrow1$ in \eqref{eq:pf-special},
then we obtain the following corollary:
\begin{corollary}
Let $n\geq1$ and $r$ be integers.
Then we have
\begin{align}
&\Pf\left((j-i)\frac{(\alpha+1)_{i+j+r-2}}{(\alpha+\beta+2)_{i+j+r-2}}\right)_{1\leq i,j\leq 2n}
\nonumber\\&\qquad=
\prod_{k=1}^{n-1}(\beta+1)_{2k}
\prod_{k=1}^n\frac{(2k-1)!(\alpha+1)_{2k+r-1}}
{(\alpha+\beta+2)_{2(k+n)+r-3}},
\label{eq:rf-ver}
\end{align}
where we use the notation
\[
(\alpha)_{n}=\begin{cases}
\prod_{i=1}^{n}(\alpha+i-1)
&\text{ if $n\geq0$,}\\
1/\prod_{i=1}^{-n}(\alpha+i+n-1)
&\text{ if $n<0$.}
\end{cases}
\]
\end{corollary}
An almost equivalent result is obtained in \cite[Theorem~6]{K1},
which is motivated by work in \cite{Ba,MW}.
In \cite{CK} Ciucu and Krattenthaler use a special case of this Pfaffian
for application to certain exact enumeration of lozenge tiling.
Further, if we put $\alpha=-\frac12$ and $\beta=\frac12$ in \eqref{eq:rf-ver},
then we obtain
\begin{align}
&\Pf\biggl((j-i)C_{i+j+r-2}\biggr)_{1\leq i,j\leq 2n}
\nonumber\\&\qquad=
\prod_{k=1}^{n-1}\frac{(4k+1)!}{(2k)!}
\prod_{k=1}^n\frac{(2k-1)!(4k+2r-2)!}
{(2k+r-1)!\{2(k+n)+r-2\}!},
\label{eq:Catalan}
\end{align}
where $C_{n}=\frac1{n+1}\binom{2n}{n}$ denotes the \defterm{Catalan numbers}.
On the other hand,
if we put $\alpha=-\frac12$ and $\beta=-\frac12$ in \eqref{eq:rf-ver},
then we obtain
\begin{align}
&\Pf\biggl((j-i)C^{(D)}_{i+j+r-2}\biggr)_{1\leq i,j\leq 2n}
\nonumber\\&\qquad=
\prod_{k=1}^{n-1}\frac{(4k)!}{(2k)!}
\prod_{k=1}^n\frac{(2k-1)!(4k+2r-2)!}
{(2k+r-1)!\{2(k+n)+r-3\}!},
\end{align}
where $C^{(D)}_{n}=\binom{2n}{n}$ is usually called the \defterm{central binomial coefficients}.

The {\sl Laguerre polynomials} (see \cite{KLS}) are defined by
\begin{equation*}
L^{(\alpha)}_{n}(x)
=\frac{(\alpha+1)_{n}}{n!}\,{}_{1}F_{1}\biggl(
{{-n}\atop{\alpha+1}}\,;\,x\biggr),
\end{equation*}
which are orthogonal with respect to the inner product
\begin{equation*}
\langle f,g\rangle=\frac1{\Gamma(\alpha+1)}\int_{0}^{\infty}e^{-x}x^{\alpha}f(x)g(x)\,dx.
\end{equation*}
Note that 
\begin{equation}
\Pf(c_{i}c_{j}a_{j}^{i})_{1\leq i,j\leq 2n}=c_{1}\ldots c_{2n}\Pf(a_{j}^{i})_{1\leq i,j\leq 2n}.
\label{eq:rc-mul}
\end{equation}
Multiplying \eqref{eq:rf-ver} by $\beta^{n(2n+1)+n(r-2)}$ and then letting $\beta\to \infty$ we get the following result.
\begin{corollary}
Let $\mu_{n}=(\alpha+1)_{n}$ for $n\geq0$,
which is known to be the moment sequence of Laguerre  polynomials.
Then we have
\begin{align}
\Pf\Bigl(
(j-i)\mu_{i+j+r-2}
\Bigr)_{1\leq i,j\leq 2n}
=\prod_{k=1}^{n}(2k-1)!(\alpha+1)_{2k+r-1}.
\label{eq:lag}
\end{align}
\end{corollary}

The {\sl Hermite polynomials} (see \cite{KLS}) are defined by
\begin{equation*}
H_{n}(x)
=(2x)^n\,{}_{2}F_{0}\biggl(
{{-n/2,-(n-1)/2}\atop{-}}\,;\,-\frac1{x^2}\biggr),
\end{equation*}
which are orthogonal with respect to the inner product
\begin{equation*}
\langle f,g\rangle=\frac1{\sqrt{\pi}}\int_{-\infty}^{\infty}e^{-x^2}f(x)g(x)\,dx.
\end{equation*}
Substituting $\alpha=\frac{1}{2}$ in \eqref{eq:lag} we get  another remarkable formula.
\begin{corollary}
Let $\mu_{n}=\prod_{k=0}^{n}(2k+1)$ denote the double factorial of $2n+1$ for $n\geq0$,
which is known to be the moment sequence of Hermite polynomials.
Then we have
\begin{align}
\Pf\Bigl(
(j-i)\mu_{i+j+r-2}
\Bigr)_{1\leq i,j\leq 2n}
=\frac1{2^n}\prod_{k=1}^{n}(4k-2)!!(4k+2r-1)!!.
\end{align}
\end{corollary}
%
%
%
%
%

%% file: tagawa04.tex
%
%%%% ---------------< ooo >--------------- %%%%
%%%%
%%%% Section 4
%%%%
%%%% ---------------< ooo >--------------- %%%%
%
%
%
\section{Proof of Theorem~\ref{conj:decomposition}}
\label{sec:proof}
%
%
%Assume $i\leq j$.
Let $a^{i}_{j}$, $t_i$ and $v^{i}_{j}$ be as in Theorem~\ref{conj:decomposition}.
To prove \eqref{eq:decomposition},
it is enough to show %that
\begin{align}
\sum_{k\geq1}\left(
v_{i}^{2k-1}t_{2k-1}v_{j}^{2k}
-v_{i}^{2k}t_{2k-1}v_{j}^{2k-1}
\right)
=a^{i}_{j}
\label{eq:plucker0}
\end{align}
for $i,j\geq1$.
%
%
%
%% ** %% which is equivalent to
%% ** %% \begin{align*}
%% ** %% &\sum_{k\geq1}\biggl\{
%% ** %% a^{2(k-1)}
%% ** %% q^{2(k-1)(2k+r-1)+1}
%% ** %% \cdot
%% ** %% \frac{
%% ** %% (q^{l-2k+1};q)_{2k-1}
%% ** %% (q^{m-2k};q)_{1}
%% ** %% (q^{m-2k+2};q)_{2k-2}
%% ** %% }{
%% ** %% (q;q)_{2k-1}
%% ** %% }
%% ** %% \\&
%% ** %% \cdot
%% ** %% \frac{
%% ** %% (aq^{2k+r};q)_{l-2k}
%% ** %% (aq;q)_{m+r-1}
%% ** %% (bq;q)_{2(k-1)}
%% ** %% }
%% ** %% {
%% ** %% (abq^{2k+r-1};q)_{2k-1}
%% ** %% (abq^{4k+r-1};q)_{l-2k}
%% ** %% (abq^{2};q)_{m+2k+r-2}
%% ** %% }
%% ** %% \cdot
%% ** %% f(2k,m,r)
%% ** %% \\&
%%%%%%%%%%%%%%%%%%%%%%%%%%%%%%%%%%%
%% ** %% -a^{2(k-1)}
%% ** %% q^{2(k-1)(2k+r-1)+1}
%% ** %% \frac{
%% ** %% (q^{l-2k};q)_{1}
%% ** %% (q^{l-2k+2};q)_{2k-2}
%% ** %% (q^{m-2k+1};q)_{2k-1}
%% ** %% }{
%% ** %% (q;q)_{2k-1}
%% ** %% }
%% ** %% \\&
%% ** %% \cdot\frac{
%% ** %% (aq^{2k+r};q)_{l-2k}
%% ** %% (aq;q)_{m+r-1}
%% ** %% (bq;q)_{2(k-1)}
%% ** %% }
%% ** %% {
%% ** %% (abq^{2k+r-1};q)_{2k-1}
%% ** %% (abq^{4k+r-1};q)_{l-2k+1}
%% ** %% (abq^{2};q)_{m+2k+r-3}
%% ** %% }
%% ** %% \cdot 
%% ** %% f(2k,l,r)
%% ** %% \biggr\}
%% ** %% =a^{l}_{m}
%% ** %% \end{align*}
%% ** %% for $l,m\geq1$.
%
%
%
%% ** %% Hence we have to show that
Replacing $aq^{r}$ by $a$,
we may assume $r=0$ hereafter without loss of generality.
Hence \eqref{eq:plucker0} is written as
\begin{align*}
&\sum_{k\geq1}
a^{2(k-1)}q^{2(k-1)(2k-1)+1}
\cdot
\frac{(q^{i-2k+2};q)_{2k-2}(q^{j-2k+2};q)_{2k-2}}{(q;q)_{2k-1}}
\nonumber\\&\cdot
\frac{
(aq^{2k};q)_{i-2k}(aq;q)_{j-1}(bq;q)_{2(k-1)}
}
{
(abq^{2k-1};q)_{2k-1}
(abq^{4k-1};q)_{i-2k+1}
(abq^2;q)_{j+2k-2}
}
\nonumber\\&\cdot
\Bigl\{
(1-q^{i-2k+1})(1-q^{j-2k})(1-abq^{i+2k-1})f(2k,j,0)
\nonumber\\&\qquad
-(1-q^{i-2k})(1-q^{j-2k+1})(1-abq^{j+2k-1})f(2k,i,0)
\Bigr\}
%\nonumber\\&
=a^{i}_{j}.
\end{align*}
Replacing $2k-1$ by $k$,
we obtain
\begin{align}
&\sum_{{k\geq1}\atop{k\text{ odd}}}
a^{k-1}q^{k(k-1)+1}
\cdot
\frac{
(abq^{2k};q)_{1}
(abq^2;q)_{k-2}
(bq,q^{i-k+1},q^{j-k+1};q)_{k-1}
}{
(q,aq,abq^{i+1},abq^{j+1};q)_{k}
}
\nonumber\\&\qquad\times
g_{k}(i,j;a,b,q)
%%\nonumber\\&
=\frac{
%%(q^{i-1}-q^{j-1})
(aq;q)_{i+j-2}
(abq^2;q)_{i-1}
(abq^2;q)_{j-1}
}{
(aq;q)_{i-1}
(aq;q)_{j-1}
(abq^2;q)_{i+j-2}
},
\label{eq:k:odd}
\end{align}
where $g_{k}(i,j;a,b,q)$ is set to be
%% 07/11/2010 %% \begin{align}
%% 07/11/2010 %% &g_{k}(i,j;a,b,q)=
%% 07/11/2010 %% (1-q^{k})(1-aq^{k})\Bigl\{
%% 07/11/2010 %% \nonumber\\&\qquad
%% 07/11/2010 %% (1-q^{i-k})(1-q^{j-k-1})
%% 07/11/2010 %% %\\&\qquad\cdot
%% 07/11/2010 %% (1-abq^{i+k})(1-abq^{j+k-1})
%% 07/11/2010 %% \nonumber\\&\qquad\qquad
%% 07/11/2010 %% -(1-q^{i-k-1})(1-q^{j-k})
%% 07/11/2010 %% %\\&\qquad\cdot
%% 07/11/2010 %% (1-abq^{i+k-1})(1-abq^{j+k})\Bigr\}/(1-q)
%% 07/11/2010 %% \nonumber\\&\qquad
%% 07/11/2010 %% +aq^{k-2}(1-b)\left(q^{i}-q^{j}\right)
%% 07/11/2010 %% \left(1-q^{i-k}\right)\left(1-q^{j-k}\right)
%% 07/11/2010 %% \left(1-abq^{2k+1}\right).
%% 07/11/2010 %% \label{eq:g_{k}(i,j;a,b,q)}
%% 07/11/2010 %% \end{align}
\begin{align}
&g_{k}(i,j;a,b,q)=(1-q^{k})(1-aq^{k})
\nonumber\\&\qquad\times
\Bigl\{
q^{-k}(1+abq^{2k})(1+abq^{i+j-1})
-ab(1+q)\left(q^{i-1}+q^{j-1}\right)\Bigr\}
\nonumber\\&\qquad
+aq^{k-1}(1-b)
\left(1-q^{i-k}\right)\left(1-q^{j-k}\right)
\left(1-abq^{2k+1}\right).
\label{eq:g_{k}(i,j;a,b,q)}
\end{align}
%%%%%%%%%%%%%%%%
By numeric experiments we observe that \eqref{eq:k:odd}
also holds in the case where the sum in the left-hand side
runs over all nonnegative even integers $k$,
%that is to say,
i.e.,
\begin{align}
&\sum_{{k\geq0}\atop{k\text{ even}}}
a^{k-1}q^{k(k-1)+1}
\cdot
\frac{
(abq^{2k};q)_{1}
(abq^2;q)_{k-2}
(bq,q^{i-k+1},q^{j-k+1};q)_{k-1}
}{
(q,aq,abq^{i+1},abq^{j+1};q)_{k}
}
\nonumber\\&\qquad\times
g_{k}(i,j;a,b,q)
%%\nonumber\\&
=\frac{
%%(q^{i-1}-q^{j-1})
(aq;q)_{i+j-2}
(abq^2;q)_{i-1}
(abq^2;q)_{j-1}
}{
(aq;q)_{i-1}
(aq;q)_{j-1}
(abq^2;q)_{i+j-2}
}.
\label{eq:k:even}
\end{align}
By adding or subtracting \eqref{eq:k:odd} and \eqref{eq:k:even},
these two identities are equivalent to
\begin{align}
&\sum_{k\geq0}
a^{k-1}q^{k(k-1)+1}
\cdot
\frac{
(abq^{2k};q)_{1}
(abq^2;q)_{k-2}
(bq,q^{i-k+1},q^{j-k+1};q)_{k-1}
}{
(q,aq,abq^{i+1},abq^{j+1};q)_{k}
}
\nonumber\\&\qquad\times
g_{k}(i,j;a,b,q)
=\frac{
2(aq;q)_{i+j-2}
(abq^2;q)_{i-1}
(abq^2;q)_{j-1}
}{
(aq;q)_{i-1}
(aq;q)_{j-1}
(abq^2;q)_{i+j-2}
},
\label{eq:add}
\end{align}
%% 07/11/2010 %% \begin{align}
%% 07/11/2010 %% &\sum_{k\geq0}
%% 07/11/2010 %% a^{k-1}q^{k(k-1)+1}
%% 07/11/2010 %% \cdot
%% 07/11/2010 %% \frac{(q^{i-k+1};q)_{k-1}(q^{j-k+1};q)_{k-1}}{(q;q)_{k}}
%% 07/11/2010 %% \nonumber\\&\qquad\times
%% 07/11/2010 %% \frac{
%% 07/11/2010 %% (bq;q)_{k-1}
%% 07/11/2010 %% (abq^2;q)_{k-2}
%% 07/11/2010 %% (abq^{2k};q)_{1}
%% 07/11/2010 %% }
%% 07/11/2010 %% {
%% 07/11/2010 %% (aq;q)_{k}
%% 07/11/2010 %% (abq^{i+1};q)_{k}
%% 07/11/2010 %% (abq^{j+1};q)_{k}
%% 07/11/2010 %% }
%% 07/11/2010 %% %%\\&
%% 07/11/2010 %% \cdot
%% 07/11/2010 %% g_{k}(i,j;a,b,q)
%% 07/11/2010 %% \nonumber\\&
%% 07/11/2010 %% =\frac{2(q^{i-1}-q^{j-1})
%% 07/11/2010 %% (aq;q)_{i+j-2}
%% 07/11/2010 %% (abq^2;q)_{i-1}
%% 07/11/2010 %% (abq^2;q)_{j-1}
%% 07/11/2010 %% }{
%% 07/11/2010 %% (aq;q)_{i-1}
%% 07/11/2010 %% (aq;q)_{j-1}
%% 07/11/2010 %% (abq^2;q)_{i+j-2}
%% 07/11/2010 %% },
%% 07/11/2010 %% \label{eq:add}
%% 07/11/2010 %% \end{align}
and
\begin{align}
&\sum_{k\geq0}
(-1)^ka^{k-1}q^{k(k-1)+1}
\cdot
\frac{
(abq^{2k};q)_{1}
(abq^2;q)_{k-2}
(bq,q^{i-k+1},q^{j-k+1};q)_{k-1}
}{
(q,aq,abq^{i+1},abq^{j+1};q)_{k}
}
\nonumber\\&\qquad\qquad\qquad\qquad\qquad\times
g_{k}(i,j;a,b,q)
=0.
\label{eq:subtract}
\end{align}
%% 07/11/2010 %% \begin{align}
%% 07/11/2010 %% &\sum_{k\geq0}
%% 07/11/2010 %% (-1)^ka^{k-1}q^{k(k-1)+1}
%% 07/11/2010 %% \cdot
%% 07/11/2010 %% \frac{(q^{i-k+1};q)_{k-1}(q^{j-k+1};q)_{k-1}}{(q;q)_{k}}
%% 07/11/2010 %% \nonumber\\&\qquad\times
%% 07/11/2010 %% \frac{
%% 07/11/2010 %% (bq;q)_{k-1}
%% 07/11/2010 %% (abq^2;q)_{k-2}
%% 07/11/2010 %% (abq^{2k};q)_{1}
%% 07/11/2010 %% }
%% 07/11/2010 %% {
%% 07/11/2010 %% (aq;q)_{k}
%% 07/11/2010 %% (abq^{i+1};q)_{k}
%% 07/11/2010 %% (abq^{j+1};q)_{k}
%% 07/11/2010 %% }
%% 07/11/2010 %% %%\\&
%% 07/11/2010 %% \cdot
%% 07/11/2010 %% g_{k}(i,j;a,b,q)
%% 07/11/2010 %% %%\nonumber\\&
%% 07/11/2010 %% =0.
%% 07/11/2010 %% \label{eq:subtract}
%% 07/11/2010 %% \end{align}
To prove \eqref{eq:add}
we rewrite $g_{k}(i,j;a,b,q)$ as follows and apply $q$-Dougall formula, i.e.,
Lemma~\ref{prop:sum_formula}, to each term,
then a direct computation leads to the desired identity:
\begin{align*}
&g_{k}(i,j;a,b,q)
=
q^{-1-k}(q+aq^{i+j})(1-q^k)(1-q^{k-1})(1-abq^k)(1-abq^{k+1})\\
&\quad
+q^{-1}
\Bigl\{
a(bq-ab-1+b)(1-q^{i+j})+(1-a)(q-ab)\\
&\qquad
+a(1+bq)(1-q^i)(1-q^j)+(1+aq^{i+j-1})(1-q)(1-abq)
\Bigr\}(1-q^k)(1-abq^k)\\
&\quad
+aq^{k-1}(1-b)(1-abq)(1-q^i)(1-q^j).
\end{align*}
To prove \eqref{eq:subtract},
we generalize this identity as
\begin{align}
&\sum_{k=0}^m
(-1)^k a^{k-1}q^{k(k-1)+1}
\nonumber\\&\qquad\times
\frac{
(1-abq^{2k})
(abq^2;q)_{k-2}
(bq,cq^{-k+1},dq^{-k+1};q)_{k-1}
\widehat g_k(a,b,c,d,q)}
{(q,aq,abcq,abdq;q)_k}
\nonumber\\&=
\frac{a^m c^{m} d^{m}(1-abq^{2m+1})
(abq^2;q)_{m-1}(bq,q/c,q/d\,;q)_{m}}
{(-q)^{m}(q,aq,abcq,abdq\,;q)_m}.
\label{eq:qv4}
\end{align}
where
\begin{align}
&\widehat g_{k}(a,b,c,d,q)
=(1-q^k)(1-aq^k)
\nonumber\\&\qquad\times
\biggl\{q^{-k}(1+abq^{2k})(1+abcdq^{-1})-abq^{-1}(1+q)(c+d)\biggr\}
\nonumber\\&\qquad
+aq^{k-1}(1-b)(1-cq^{-k})(1-dq^{-k})(1-abq^{2k+1}).
\end{align}
Then \eqref{eq:qv4} is proven by induction on $m$.
This completes the proof of Theorem~\ref{conj:decomposition}.
%
%
%
%
%% ------------ < *** > ------------ %%
%% Lemma
%% ------------ < *** > ------------ %%
\begin{lemma}
\label{prop:sum_formula}
Let $m$ be an integer.
Then we have
%% 05/11/2010%% \begin{align}
%% 05/11/2010%% &\sum_{k\geq m}
%% 05/11/2010%% a^{k-m}q^{k(k-m)}
%% 05/11/2010%% \cdot
%% 05/11/2010%% \frac{(q^{i-k+1};q)_{k-1}(q^{j-k+1};q)_{k-1}}{(q;q)_{k-m}}
%% 05/11/2010%% \nonumber\\&\qquad\qquad\times
%% 05/11/2010%% \frac{
%% 05/11/2010%% (bq;q)_{k-1}
%% 05/11/2010%% (abq^2;q)_{k+m-2}
%% 05/11/2010%% (1-abq^{2k})
%% 05/11/2010%% }
%% 05/11/2010%% {
%% 05/11/2010%% (aq;q)_{k}
%% 05/11/2010%% (abq^{i+1};q)_{k}
%% 05/11/2010%% (abq^{j+1};q)_{k}
%% 05/11/2010%% }
%% 05/11/2010%% %%\\&
%% 05/11/2010%% \nonumber\\&
%% 05/11/2010%% =(q^{i-m+1},q^{j-m+1},bq;q)_{m-1}\cdot
%% 05/11/2010%% \frac{
%% 05/11/2010%% (aq;q)_{i+j-m}
%% 05/11/2010%% (abq^2;q)_{i-1}
%% 05/11/2010%% (abq^2;q)_{j-1}
%% 05/11/2010%% }{
%% 05/11/2010%% (aq;q)_{i}
%% 05/11/2010%% (aq;q)_{j}
%% 05/11/2010%% (abq^2;q)_{i+j-1}
%% 05/11/2010%% }.
%% 05/11/2010%% \label{eq:Dougall}
%% 05/11/2010%% \end{align}
\begin{align}
&\sum_{k\geq m}
a^{k-m}q^{k(k-m)}
\cdot
\frac{
(1-abq^{2k})
(abq^2;q)_{k+m-2}
(bq,q^{i-k+1},q^{j-k+1};q)_{k-1}
}
{
(q;q)_{k-m}
(aq,abq^{i+1},abq^{j+1};q)_{k}
}
\nonumber\\&
=(q^{i-m+1},q^{j-m+1},bq;q)_{m-1}\cdot
\frac{
(aq^{j+1};q)_{i-m}
(abq^2;q)_{i-1}
}{
(aq;q)_{i}
(abq^{j+1};q)_{i}
}.
\label{eq:Dougall}
\end{align}
\end{lemma}
%
%
%
%% ------------ < *** > ------------ %%
%% Dougall's formula and Bailey pair
%% ------------ < *** > ------------ %%
In fact \eqref{eq:Dougall} reduces to the $q$-Dougall formula (Jackson's formula)
\cite[(12.3.2)]{AAR}, \cite[(2.4.2)]{GR}
\begin{align}
{}_{6}\phi_{5}\left[{
{a,qa^{\frac12},-qa^{\frac12},b,c,q^{-n}}
\atop
{a^{\frac12},-a^{\frac12},aq/b,aq/c,aq^{n+1}}
}\,;\,q,\frac{aq^{n+1}}{bc}
\right]
=\frac{(aq,aq/bc;q)_{n}}{(aq/b,aq/c;q)_{n}},
\label{eq:q-Dougall}
\end{align}
by the substitution
\[
a\leftarrow abq^{2m},\quad
b\leftarrow bq^{m},\quad
c\leftarrow q^{m-j},\quad
n\leftarrow i-m.
\]
\begin{remark}
The lemma also directly follows from the Bailey pair $(\alpha_{n},\beta_{n})$
given by
\begin{align*}
&\alpha_{n}
=\frac{(a,b,c;q)_{n}(1-aq^{2n})(aq/bc)^n(-1)^nq^{\binom{n}2}}
{(q,aq/b,aq/c;q)_{n}(1-a)},
\\
&\beta_{n}
=\frac
{(aq/bc;q)_{n}}
{(q,aq/b,aq/c;q)_{n}}.
\end{align*}
Here a pair $(\alpha_{n},\beta_{n})$ is said to be a \defterm{Bailey pair} \cite{AAR} if it satisfies
\[
\beta_{n}=\sum_{k=0}^{n}\frac{\alpha_{k}}{(q;q)_{n-k}(aq;q)_{n+k}}.
\]
\end{remark}
%%\medbreak

In  fact we prove two identities \eqref{eq:k:odd} and \eqref{eq:k:even} in this section.
While \eqref{eq:k:odd} is used to prove Theorem~\ref{conj:decomposition},
one may ask what's \eqref{eq:k:even} for?
In fact we can interpret \eqref{eq:k:even} as a Pfaffian decomposition of another skew-symmetric matrix.
Define $a^{i}_{j}$ for $i,j\geq0$ by
\begin{equation}
a^{0}_{j}=\frac{(abq^{r-1};q)_{1}(aq;q)_{j+r-1}}{aq^{r}(1-b)(abq^{2};q)_{j+r-2}}
\label{eq:a^{0}_{j}},
\end{equation}
with $a^{i}_{0}=-a^{0}_{i}$, $a^{0}_{0}=0$
and \eqref{eq:a^{i}_{j}} for $i,j\geq1$.

\begin{theorem}
\label{th:byproduct}
Let $\widetilde A=(a^{i}_{j})_{i,j\geq0}$ where $a^{i}_{j}$ is as above.
Let $t_{i}$, $o^{i}_{j}$ and $e^{i}_{j}$ be as in Theorem~\ref{conj:decomposition},
and we put
\begin{align*}
{\widetilde v}^{i}_{j}
&=\begin{cases}
o^{i}_{j}&\text{ if $i$ is even,}\\
e^{i}_{j}&\text{ if $i$ is odd.}
\end{cases}
\end{align*}
If we set
$\displaystyle
\widetilde T=\bigoplus_{i\geq0}\begin{pmatrix} 0&t_{2i}\\-t_{2i}&0\end{pmatrix}
$
and
$
\widetilde V=\left({\widetilde v}^{i}_{j}\right)_{i,j\geq0}
$
then
\begin{equation}
\widetilde A={}^{t}{\widetilde V}\,\widetilde T\, {\widetilde V}
\label{eq:by-decomposition}
\end{equation}
gives the Pfaffian decomposition of $\widetilde A$.
\end{theorem}
%
%
%
%% ------------------------- %%
%% Catalan Pfaffian
%% ------------------------- %%
%
%
\begin{corollary}
\label{cor:byproduct}
Let $n\geq1$ and $r$ be integers.
Then we have
\begin{align}
\Pf\left(a^{i}_{j}\right)_{0\leq i,j\leq 2n-1}
&=a^{n(n-2)}q^{n(n-1)(4n-5)/3+n(n-2)r}
\nonumber\\&\times
\prod_{k=0}^{n-1}
\frac
{(q;q)_{2k}(aq;q)_{2k+r}(bq;q)_{2k-1}}
{(abq^2;q)_{4k+r-1}(abq^{2k+r};q)_{2k-1}}.
\label{eq:pf-byproduct}
\end{align}
Let $P_{n,r}(a,b;q)$ denote the right-hand side of \eqref{eq:pf-byproduct}.
Then, more generally we have
\begin{align}
\Pf\left({\widetilde A}_{[0,2n-2],m-1}\right)
&=
%\nonumber\\&\times
\frac{(q^{m-2n+1};q)_{2n-2}(aq^{2n+r-1};q)_{m-2n}}
{(q;q)_{2n-2}(abq^{4n+r-3};q)_{m-2n}}
P_{n,r}(a,b;q),
\label{eq:pf-general3}
\\
\Pf\left({\widetilde A}_{[0,2n-3],2n-1,m-1}\right)
&=
q\cdot\frac{
(q^{m-2n};q)_{1}
(q^{m-2n+2};q)_{2n-3}
(aq^{2n+r-1};q)_{m-2n}
}
{
(q;q)_{2n-2}
(abq^{4n+r-5};q)_{1}
(abq^{4n+r-3};q)_{m-2n+1}
}
\nonumber\\&\times
f(2n-1,m-1,r)
P_{n,r}(a,b;q)
.
\label{eq:pf-general4}
\end{align}
Here we use the notation $[i,j]=\{\,x\in\mathbb{Z}\,|\,i\leq x\leq j\,\}$.
\end{corollary}
%

%% file: tagawa05.tex
%
%%%% ---------------< ooo >--------------- %%%%
%%%%
%%%% Section 5
%%%%
%%%% ---------------< ooo >--------------- %%%%
%
%
%
\section{Weighted enumeration of shifted RPPs}
\label{sec:application}
%
%
%Assume $i\leq j$.
In this section we give an application of Corollary~\ref{conj:01},
which enumerates a certain class of shifted reverse plane partitions.
%
%
%
%% ------------------------- %%
%% Shifted Tableaux
%% ------------------------- %%
%
%
\begin{definition}
\label{def:SRPP}
A \defterm{shifted reverse plane partition} (abbreviated as \defterm{shifted RPP}) is an array 
$
\pi=(\pi_{ij})
$
of nonnegative integers,
defined only for $i\leq j$,
that has nondecreasing rows and columns,
and that can be written in the form
\begin{equation}
\begin{array}{ccccccc}
\pi_{11}&\pi_{12}&\pi_{13}&\hdots&\hdots&\hdots&\pi_{1,\lambda_{1}}\\
        &\pi_{22}&\pi_{23}&\hdots&\hdots&\pi_{2,\lambda_{2}+1}\\
        &        &\ddots&\vdots&\vdots&\rdots\\
        &        &      &\pi_{n,n}&\hdots&\pi_{n,\lambda_{n}+n-1}\\
\end{array},
\label{eq:ShiftedTableau}
\end{equation}
where
\begin{enumerate}
\item[(i)]
$\lambda_{1}>\lambda_{2}>\dots>\lambda_{n}>0$,
\item[(ii)]
$\pi_{i,j}\leq\pi_{i,j+1}$ whenever the both sides are defined,
\item[(iii)]
$\pi_{i,j}\leq\pi_{i+1,j}$ whenever the both sides are defined.
\end{enumerate}
Further, 
if $\pi$ also satisfies
\begin{enumerate}
\item[(iii')]
$\pi_{i,j}<\pi_{i+1,j}$ whenever the both sides are defined,
\end{enumerate}
then it is called \defterm{column-strict}
shifted reverse plane partition or a \defterm{shifted tableau}.
The entries $\pi_{ij}$ are called the \defterm{parts} of $\pi$.
To each shifted reverse plane partition $\pi$
we assign the \defterm{weight} $|\pi|=\sum_{ij}\pi_{ij}$ to be the sum of parts.
The strict partition $\lambda$ is called the \defterm{shape} of $\pi$,
and the nondecreasing sequence $(\pi_{11},\pi_{22},\dots,\pi_{nn})$
is called the \defterm{profile} of $\pi$.
Let $\mathcal{R}_{\lambda,\mu}$ denote the set of all shifted reverse plane partitions
of shape $\lambda$ and profile $\mu$,
and $\mathcal{T}_{\lambda,\mu}$ the set of all shifted tableaux
of shape $\lambda$ and profile $\mu$
for fixed $\lambda=(\lambda_{1},\dots,\lambda_{n})$ and 
$\mu=(\mu_{1},\dots,\mu_{n})$ with 
$\lambda_{1}>\dots>\lambda_{n}>0$
and $0\leq\mu_{1}\leq\dots\leq\mu_{n}$.
\end{definition}
%
%
%
%% ------------------------------------- %%
%% 
%% ------------------------------------- %%
For example,
\[
\smyoung{
0&0&0&0&1&1&1&2&2&3&4\\
\blank&1&1&2&2&2&3&3&4\\
\blank&\blank&3&3&3&4&4&5\\
\blank&\blank&\blank&5&5&5\\
}
\]
is a shifted tableau of shape $\lambda=(11,8,6,3)$ and profile $\mu=(0,1,3,5)$
with weight $69$.
%\smallskip
%
If $\mathcal{F}$ is a family of shifted reverse plane partitions,
then the generating function of $\mathcal{F}$ is defined to be
\begin{equation}
\GF{\mathcal{F}}
=\sum_{\pi\in\mathcal{F}}q^{|\pi|}.
\label{eq:GeneratingFunction}
\end{equation}
%
%
%% ------------------------------------- %%
%% RPP and tableau
%% ------------------------------------- %%
Let $n(\lambda)=\sum_{i\geq1}(i-1)\lambda_{i}$.
It is easy to see that
\begin{equation}
\GF{\mathcal{T}_{\lambda,\nu+\epsilon_{n}}}
=q^{n(\lambda)}\GF{\mathcal{\mathcal{R}_{\lambda,\nu}}},
\label{eq:GF-rel}
\end{equation}
where $\epsilon_{n}=(0,1,\dots,n-1)$ and $\nu=(\nu_{1},\nu_{2},\dots,\nu_{n})$
is a profile such that $\nu_{1}\leq\nu_{2}\leq\dots\leq\nu_{n}$.
%
%
%% ------------------------------------- %%
%% Profile and weight
%% ------------------------------------- %%
Let $\mathcal{P}_{n}$ denote the set of profiles
$\nu=(\nu_{1},\dots,\nu_{2n})$ such that $0\leq\nu_{1}\leq\nu_{2}\leq\dots\leq\nu_{2n}$
and $\nu_{2k}=\nu_{2k-1}$ for $k=1,\dots,n$.
For $\nu\in\mathcal{P}_{n}$ and $x\in\mathbb{Z}$ we let
\begin{equation}
\omega_{x}(\nu)
=\left(aq^{x}\right)^{|\nu|/2}
\prod_{k=1}^{n}\frac{(bq^{2k-1};q)_{\nu_{2k-1}}}{(q^{2k-1};q)_{\nu_{2k-1}}},
\label{eq:profile-weight}
\end{equation}
where $|\nu|=\sum_{k=1}^{2n}\nu_{k}$.
Let $\mathcal{P}_{n}'=\left\{\nu+\epsilon_{2n}\,|\,\nu\in\mathcal{P}_{n}\right\}$,
and
\begin{equation}
\omega_{x}'(\mu)
=\left(aq^{x}\right)^{\left(|\mu|-n\right)/2}
\prod_{k=1}^{n}\frac{(bq;q)_{\mu_{2k-1}}}{(q;q)_{\mu_{2k-1}}}
\end{equation}
for $\mu\in\mathcal{P}_{n}'$.
%% ------------------------------------- %%
%% 
%% ------------------------------------- %%
Now we are in position to state our main theorem in this section.
If the shape $\lambda=(\lambda_{1},\lambda_{2},\dots,\lambda_{r})$ is in the form of
$\lambda=(m,m-1,\dots,m-r+1)$ for positive integers $m\geq r$,
then it is called \defterm{staircase}.
\begin{figure}[htb]
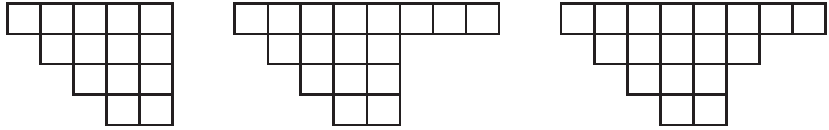

\[
\smyoung{
&&&&\\
\blank&&&&\\
\blank&\blank&&&\\
\blank&\blank&\blank&&\\
}
\qquad
\smyoung{
&&&&&&&\\
\blank&&&&\\
\blank&\blank&&&\\
\blank&\blank&\blank&&\\
}
\qquad
\smyoung{
&&&&&&&\\
\blank&&&&&\\
\blank&\blank&&&\\
\blank&\blank&\blank&&\\
}
\]
\caption{Nearly Staircase Shapes}\label{fig:StaircaseShape}
\end{figure}
Each of \eqref{eq:pf-special}, \eqref{eq:pf-general1} and \eqref{eq:pf-general2} 
corresponds to each of \eqref{eq:ShiftedRPP-GF}, \eqref{eq:ShiftedRPP-GF2} and \eqref{eq:ShiftedRPP-GF3} below.
In fact the leftmost diagram in Figure~\ref{fig:StaircaseShape}
gives the case of $m=5$, $n=2$ in \eqref{eq:ShiftedRPP-GF},
and the middle (resp. rightmost) diagram in Figure~\ref{fig:StaircaseShape}
gives the case of $l=8$, $m=5$, $n=2$ in \eqref{eq:ShiftedRPP-GF2}
(resp. \eqref{eq:ShiftedRPP-GF3}).
%
%
%% -------------------------------------------------- %%
%% Theorem (Weighted enumration of shifted tableaux)
%% -------------------------------------------------- %%
%
%
\begin{theorem}
\label{th:ShiftedTableaux}
Let $r$ be an integer.
For any positive integers $m$ and $n$ such that $m\geq2n$,
we fix the shape $\lambda=(m,m-1,\dots,m-2n+1)$ of length $2n$.
Then we have
\begin{align}
&\sum_{\nu\in\mathcal{P}_{n}}
\omega_{r-2(m-2n)-1}(\nu)\,
\GF{\mathcal{R}_{\lambda,\nu}}
\nonumber\\&=
\left\{
\frac{(abq^2;q)_{\infty}}{(aq;q)_{\infty}}
\right\}^{n}
\cdot
\prod_{k=1}^{2n}
\frac{(q;q)_{k-1}}{(q;q)_{k+m-2n-1}}
\cdot
\prod_{k=1}^{n}
\frac{(aq;q)_{2k+r-1}}
{(abq^2;q)_{2(k+n)+r-3}}.
\label{eq:ShiftedRPP-GF}
\end{align}
More generally,
if $\lambda=(l,m-1,m-2,m-3,\dots,m-2n+1)$ where $l\geq m$,
then we have
\begin{align}
&\sum_{\nu\in\mathcal{P}_{n}}
\omega_{r-2(m-2n)-1}(\nu)\,
\GF{\mathcal{R}_{\lambda,\nu}}
=\left\{
\frac{(abq^2;q)_{\infty}}{(aq;q)_{\infty}}
\right\}^{n}
\cdot
\frac{(q^{l-m+1};q)_{2n-1}}
{(q;q)_{l-1}}
\nonumber\\&\times
\prod_{k=1}^{2n-1}\frac{(q;q)_{k-1}}{(q;q)_{k+m-2n-1}}
\cdot
\frac{(aq;q)_{l-m+2n+r-1}}
{(abq^2;q)_{l-m+4n+r-3}}
\cdot
\prod_{k=1}^{n-1}\frac{
(aq;q)_{2k+r-1}
}
{(abq^2;q)_{2(k+n)+r-3}},
\label{eq:ShiftedRPP-GF2}
\end{align}
and
if $\lambda=(l,m,m-2,m-3,\dots,m-2n+1)$ where $l>m$,
 then we have
\begin{align}
&\sum_{\nu\in\mathcal{P}_{n}}
\omega_{r-2(m-2n)-1}(\nu)\,
\GF{\mathcal{R}_{\lambda,\nu}}
\nonumber\\&
=\left\{
\frac{(abq^2;q)_{\infty}}{(aq;q)_{\infty}}
\right\}^{n}
\cdot
\frac{(q^{l-m};q)_{1}(q^{l-m+2};q)_{2n-2}}{(q;q)_{l-1}(q;q)_{m-1}}
\cdot
\frac{\prod_{k=1}^{2n-1}(q;q)_{k-1}}{\prod_{k=1}^{2n-2}(q;q)_{k+m-2n-1}}
\nonumber\\&\times
\frac{(aq;q)_{l-m+2n+r-1}f(2n,l-m+2n,r)}
{(abq^{4n+r-3};q)_{1}(abq^2;q)_{l-m+4n+r-2}}
\prod_{k=1}^{n-1}\frac{
(aq;q)_{2k+r-1}
}
{(abq^2;q)_{2(k+n)+r-3}}.
\label{eq:ShiftedRPP-GF3}
\end{align}
\end{theorem}
To prove this theorem,
 we first recall the notation of the lattice path method,
which is due to Gessel and Viennot \cite{GV}.
%Thus we recall the notation.
Let $D=(V,E)$ be an acyclic digraph without multiple edges.
If $u$ and $v$ are any pair of vertices,
let $\PATH{u}{v}$ denote the set of all directed paths from $u$ to $v$.
For a fixed positive integer $n$,
an \defterm{$n$-vertex} is an $n$-tuple of vertices of $D$.
If $\boldsymbol{u}=(u_1,\dots,u_n)$ and $\boldsymbol{v}=(v_1,\dots,v_n)$ are $n$-vertices,
an \defterm{$n$-path} from $\boldsymbol{u}$ to $\boldsymbol{v}$ is an $n$-tuple $\boldsymbol{P}=(P_1,\dots,P_n)$
such that $P_i\in\PATH{u_i}{v_i}$, $i=1,\dots,n$.
The $n$-path $\boldsymbol{P}=(P_1,\dots,P_n)$ is said to be \defterm{non-intersecting}
if any two different paths $P_i$ and $P_j$ have no vertex in common.
We will write $\PATH{\boldsymbol{u}}{\boldsymbol{v}}$ for the set of all $n$-paths from $\boldsymbol{u}$ to $\boldsymbol{v}$,
and write $\NPATH{\boldsymbol{u}}{\boldsymbol{v}}$ for the subset of $\PATH{\boldsymbol{u}}{\boldsymbol{v}}$ 
consisting of non-intersecting $n$-paths.
If $\boldsymbol{u}=(u_1,\dots,u_m)$ and $\boldsymbol{v}=(v_1,\dots,v_n)$ are linearly ordered sets of vertices of $D$,
then $\boldsymbol{u}$ is said to be \defterm{$D$-compatible} with $\boldsymbol{v}$ if
every path $P\in{\cal P}(u_i,v_l)$ intersects with every path $Q\in{\cal P}(u_j,v_k)$ whenever $i<j$ and $k<l$.
%Let $S_n$ denote the symmetric group on $\{1,2,\dots,n\}$.
If $\pi\in S_n$,
by $\boldsymbol{v}^\pi$ we mean the $n$-vertex $(v_{\pi(1)},\dots,v_{\pi(n)})$.
The weight $w(\boldsymbol{P})$ of an $n$-path $\boldsymbol{P}$ is defined to be the product
of the weights of its components.
Thus,
if $\boldsymbol{u}=(u_1,\dots,u_n)$ and $\boldsymbol{v}=(v_1,\dots,v_n)$ are $n$-vertices,
 we define the generating functions $F(\boldsymbol{u},\boldsymbol{v})=\GF{\PATH{\boldsymbol{u}}{\boldsymbol{v}}}=\sum_{\boldsymbol{P}\in\PATH{\boldsymbol{u}}{\boldsymbol{v}}}w(\boldsymbol{P})$ 
and $F_0(\boldsymbol{u},\boldsymbol{v})=\GF{\NPATH{\boldsymbol{u}}{\boldsymbol{v}}}=\sum_{\boldsymbol{P}\in\NPATH{\boldsymbol{u}}{\boldsymbol{v}}}w(\boldsymbol{P})$.
In particular, if $u$ and $v$ are any pair of vertices,
we write
\begin{equation*}
h(u,v)=\GF{\PATH{u}{v}}=\sum_{P\in\PATH{u}{v}}w(P).
\end{equation*}
%
%
%
%--------------------------------------------------------%
% Theorem (Lidstr\"om-Gessel-Viennot)
%--------------------------------------------------------%
\begin{lemma}
\label{lem:Gessel-Viennot}
(Lindstr\"om-Gessel-Viennot \cite{GV})

Let $\boldsymbol{u}=(u_1,\dots,u_n)$ and $\boldsymbol{v}=(v_1,\dots,v_n)$ be two $n$-vertices in an acyclic digraph $D$.
Then
\begin{equation}
\sum_{\pi\in S_n}\sgn\pi\ F_0(\boldsymbol{u}^{\pi},\boldsymbol{v})=\det[h(u_i,v_j)]_{1\le i,j\le n}.
\label{eq:LGV1}
\end{equation}
In particular, if $\boldsymbol{u}$ is $D$-compatible with $\boldsymbol{v}$,
then
\begin{equation}
F_0(\boldsymbol{u},\boldsymbol{v})=\det[h(u_i,v_j)]_{1\le i,j\le n}.
\label{eq:LGV2}
\end{equation}
\end{lemma}
Using the Lindstr\"om-Gessel-Viennot theorem, 
we obtain the following determinantal expression for the generating function
of shifted tableaux.
%% ------------------------- %%
%% Lemma for GF
%% ------------------------- %%
%
%
\begin{lemma}
\label{lem:STGF}
Let 
$\lambda=(\lambda_{1},\dots,\lambda_{n})$ and 
$\mu=(\mu_{1},\dots,\mu_{n})$ be sequences
such that $\lambda_{1}>\dots>\lambda_{n}>0$
and $0\leq\mu_{1}<\dots<\mu_{n}$.
Then
\begin{equation}
\GF{\mathcal{T}_{\lambda,\mu}}
=\det\left(
\frac{q^{\lambda_{i}\mu_{j}}}{(q;q)_{\lambda_{i}-1}}
\right)_{1\leq i,j\leq n}.
\label{eq:STGF}
\end{equation}
\end{lemma}
\begin{demo}{Proof}
We consider the digraph $D$ whose vertex set is $\mathbb{Z}_{\geq0}^2$ and the edge set is defined as follows.
An edge is directed from $u$ to $v$ whenever $v-u=(1,0)$ or $(0,1)$
(resp. whenever $v-u=(1,0)$)
if the vertex $u=(i,j)$ satisfies $i>0$ (resp. $i=0$).
For $u=(i,j)$, we assign the weight $q^j$ (resp. $1$) to the edge with $v-u=(1,0)$ (resp. $(0,1)$).
Fix a sufficiently large positive integer $N$.
For given
$\lambda=(\lambda_{1},\dots,\lambda_{n})$
 and 
$\mu=(\mu_{1},\dots,\mu_{n})$,
we set the vertices
$u_{i}=(0,\mu_{i})$ and $v_{j}=(\lambda_{j},N)$ for $i,j=1,\dots,n$.
Let $\mathcal{T}^{N}_{\lambda,\mu}$ denote the set of shifted tableaux $\pi\in\mathcal{T}_{\lambda,\mu}$
such that each part is less than or equal to $N$.
%
%
%
%% ------------< *** >------------ %%
% Figure : Lattice paths
%% ------------< *** >------------ %%
%
\begin{figure}[htb]
\begin{center}
%\hspace*{30mm}
%\setlength{\unitlength}{1pt} %default
\setlength{\unitlength}{0.5mm}
\begin{picture}(140,90)
% 
% vertical & horizontal
%
\put( 10, 10){\vector( 1, 0){130}}
\put( 10, 10){\vector( 0, 1){80}}
\put( 10, 10){\vector( 1, 0){10}}
\put( 20, 10){\vector( 1, 0){10}}
\put( 30, 10){\vector( 1, 0){10}}
\put( 40, 10){\vector( 1, 0){10}}
\put( 50, 10){\vector( 0, 1){10}}
\put( 50, 20){\vector( 1, 0){10}}
\put( 60, 20){\vector( 1, 0){10}}
\put( 70, 20){\vector( 1, 0){10}}
\put( 80, 20){\vector( 0, 1){10}}
\put( 80, 30){\vector( 1, 0){10}}
\put( 90, 30){\vector( 1, 0){10}}
\put(100, 30){\vector( 0, 1){10}}
\put(100, 40){\vector( 1, 0){10}}
\put(110, 40){\vector( 0, 1){10}}
\put(110, 50){\vector( 1, 0){10}}
\put(120, 50){\vector( 0, 1){10}}
\put(120, 60){\vector( 0, 1){10}}
\put(120, 70){\vector( 0, 1){10}}
\put( 10, 20){\vector( 1, 0){10}}
\put( 20, 20){\vector( 1, 0){10}}
\put( 30, 20){\vector( 0, 1){10}}
\put( 30, 30){\vector( 1, 0){10}}
\put( 40, 30){\vector( 1, 0){10}}
\put( 50, 30){\vector( 1, 0){10}}
\put( 60, 30){\vector( 0, 1){10}}
\put( 60, 40){\vector( 1, 0){10}}
\put( 70, 40){\vector( 1, 0){10}}
\put( 80, 40){\vector( 0, 1){10}}
\put( 80, 50){\vector( 1, 0){10}}
\put( 90, 50){\vector( 0, 1){10}}
\put( 90, 60){\vector( 0, 1){10}}
\put( 90, 70){\vector( 0, 1){10}}
\put( 10, 40){\vector( 1, 0){10}}
\put( 20, 40){\vector( 1, 0){10}}
\put( 30, 40){\vector( 1, 0){10}}
\put( 40, 40){\vector( 0, 1){10}}
\put( 40, 50){\vector( 1, 0){10}}
\put( 50, 50){\vector( 1, 0){10}}
\put( 60, 50){\vector( 0, 1){10}}
\put( 60, 60){\vector( 1, 0){10}}
\put( 70, 60){\vector( 0, 1){10}}
\put( 70, 70){\vector( 0, 1){10}}
\put( 10, 60){\vector( 1, 0){10}}
\put( 20, 60){\vector( 1, 0){10}}
\put( 30, 60){\vector( 1, 0){10}}
\put( 40, 60){\vector( 0, 1){10}}
\put( 40, 70){\vector( 0, 1){10}}
\put( 20, 10){\circle*{1.0}}
\put( 30, 10){\circle*{1.0}}
\put( 40, 10){\circle*{1.0}}
\put( 50, 10){\circle*{1.0}}
\put( 60, 10){\circle*{1.0}}
\put( 70, 10){\circle*{1.0}}
\put( 80, 10){\circle*{1.0}}
\put( 90, 10){\circle*{1.0}}
\put(100, 10){\circle*{1.0}}
\put(110, 10){\circle*{1.0}}
\put(120, 10){\circle*{1.0}}
\put(130, 10){\circle*{1.0}}
\put( 10, 20){\circle*{1.0}}
\put( 20, 20){\circle*{1.0}}
\put( 30, 20){\circle*{1.0}}
\put( 40, 20){\circle*{1.0}}
\put( 50, 20){\circle*{1.0}}
\put( 60, 20){\circle*{1.0}}
\put( 70, 20){\circle*{1.0}}
\put( 80, 20){\circle*{1.0}}
\put( 90, 20){\circle*{1.0}}
\put(100, 20){\circle*{1.0}}
\put(110, 20){\circle*{1.0}}
\put(120, 20){\circle*{1.0}}
\put(130, 20){\circle*{1.0}}
\put( 10, 30){\circle*{1.0}}
\put( 20, 30){\circle*{1.0}}
\put( 30, 30){\circle*{1.0}}
\put( 40, 30){\circle*{1.0}}
\put( 50, 30){\circle*{1.0}}
\put( 60, 30){\circle*{1.0}}
\put( 70, 30){\circle*{1.0}}
\put( 80, 30){\circle*{1.0}}
\put( 90, 30){\circle*{1.0}}
\put(100, 30){\circle*{1.0}}
\put(110, 30){\circle*{1.0}}
\put(120, 30){\circle*{1.0}}
\put(130, 30){\circle*{1.0}}
\put( 10, 40){\circle*{1.0}}
\put( 20, 40){\circle*{1.0}}
\put( 30, 40){\circle*{1.0}}
\put( 40, 40){\circle*{1.0}}
\put( 50, 40){\circle*{1.0}}
\put( 60, 40){\circle*{1.0}}
\put( 70, 40){\circle*{1.0}}
\put( 80, 40){\circle*{1.0}}
\put( 90, 40){\circle*{1.0}}
\put(100, 40){\circle*{1.0}}
\put(110, 40){\circle*{1.0}}
\put(120, 40){\circle*{1.0}}
\put(130, 40){\circle*{1.0}}
\put( 10, 50){\circle*{1.0}}
\put( 20, 50){\circle*{1.0}}
\put( 30, 50){\circle*{1.0}}
\put( 40, 50){\circle*{1.0}}
\put( 50, 50){\circle*{1.0}}
\put( 60, 50){\circle*{1.0}}
\put( 70, 50){\circle*{1.0}}
\put( 80, 50){\circle*{1.0}}
\put( 90, 50){\circle*{1.0}}
\put(100, 50){\circle*{1.0}}
\put(110, 50){\circle*{1.0}}
\put(120, 50){\circle*{1.0}}
\put(130, 50){\circle*{1.0}}
\put( 10, 60){\circle*{1.0}}
\put( 20, 60){\circle*{1.0}}
\put( 30, 60){\circle*{1.0}}
\put( 40, 60){\circle*{1.0}}
\put( 50, 60){\circle*{1.0}}
\put( 60, 60){\circle*{1.0}}
\put( 70, 60){\circle*{1.0}}
\put( 80, 60){\circle*{1.0}}
\put( 90, 60){\circle*{1.0}}
\put(100, 60){\circle*{1.0}}
\put(110, 60){\circle*{1.0}}
\put(120, 60){\circle*{1.0}}
\put(130, 60){\circle*{1.0}}
\put( 10, 70){\circle*{1.0}}
\put( 20, 70){\circle*{1.0}}
\put( 30, 70){\circle*{1.0}}
\put( 40, 70){\circle*{1.0}}
\put( 50, 70){\circle*{1.0}}
\put( 60, 70){\circle*{1.0}}
\put( 70, 70){\circle*{1.0}}
\put( 80, 70){\circle*{1.0}}
\put( 90, 70){\circle*{1.0}}
\put(100, 70){\circle*{1.0}}
\put(110, 70){\circle*{1.0}}
\put(120, 70){\circle*{1.0}}
\put(130, 70){\circle*{1.0}}
\put( 10, 80){\circle*{1.0}}
\put( 20, 80){\circle*{1.0}}
\put( 30, 80){\circle*{1.0}}
\put( 40, 80){\circle*{1.0}}
\put( 50, 80){\circle*{1.0}}
\put( 60, 80){\circle*{1.0}}
\put( 70, 80){\circle*{1.0}}
\put( 80, 80){\circle*{1.0}}
\put( 90, 80){\circle*{1.0}}
\put(100, 80){\circle*{1.0}}
\put(110, 80){\circle*{1.0}}
\put(120, 80){\circle*{1.0}}
\put(130, 80){\circle*{1.0}}
%
%  Big circles
\put( 10, 10){\circle*{2}}
\put( 10, 20){\circle*{2}}
\put( 10, 40){\circle*{2}}
\put( 10, 60){\circle*{2}}
\put( 40, 80){\circle*{2}}
\put( 70, 80){\circle*{2}}
\put( 90, 80){\circle*{2}}
\put(120, 80){\circle*{2}}
%
% labels
% vertices
\put(  0, 10){\makebox(0,0){$\scriptstyle(0,\mu_{1})$}}
\put(  0, 20){\makebox(0,0){$\scriptstyle(0,\mu_{2})$}}
\put(  0, 40){\makebox(0,0){$\scriptstyle\vdots$}}
\put(  0, 60){\makebox(0,0){$\scriptstyle(0,\mu_{r})$}}
\put( 40, 85){\makebox(0,0){$\scriptstyle(\lambda_{r},N)$}}
\put( 70, 85){\makebox(0,0){$\scriptstyle\hdots$}}
\put( 90, 85){\makebox(0,0){$\scriptstyle(\lambda_{2},N)$}}
\put(120, 85){\makebox(0,0){$\scriptstyle(\lambda_{1},N)$}}
\put(140,  5){\makebox(0,0){$\scriptstyle x$}}
\put(  5, 90){\makebox(0,0){$\scriptstyle y$}}
%
%
%%\put( 15, 18){\makebox(0,0){\footnotesize{$\scriptstyle0$}}}
%%\put( 25, 18){\makebox(0,0){\footnotesize{$\scriptstyle1$}}}
%%\put( 35, 13){\makebox(0,0){\footnotesize{$\scriptstyle0$}}}
%%\put( 45, 18){\makebox(0,0){\footnotesize{$\scriptstyle0$}}}
%
\end{picture}
\caption{The $n$-path corresponding to the above shifted tableaux}\label{fig:SPath}
\end{center}
\end{figure}
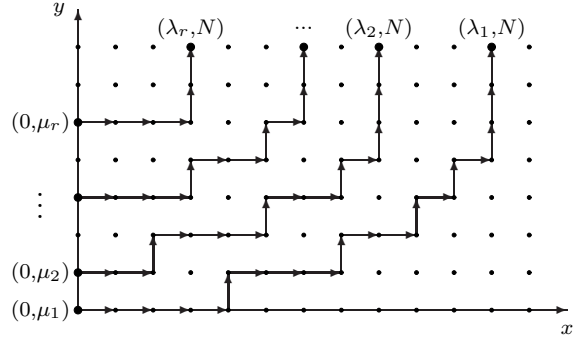
Then a tableau $\pi\in\mathcal{T}^{N}_{\lambda,\mu}$ is interpreted as
an $n$-path $\boldsymbol{P}$ from $\boldsymbol{u}=(u_{1},\dots,u_{n})$
to $\boldsymbol{v}=(v_{1},\dots,v_{n})$.
For instance the shifted tableaux in the above example is pictorially illustrated
by the $4$-path in Figure~\ref{fig:SPath}.
If $u=(0,y)$ and $v=(x,N)$,
then we have $h(u,v)=q^{xy}\qbinom{x-1+N-y}{x-1}$.
Hence, by Lemma~\ref{lem:Gessel-Viennot},
we obtain
\begin{equation*}
\GF{\mathcal{T}^{N}_{\lambda,\mu}}
=\det\left(
q^{\lambda_{i}\mu_{j}}\qbinom{\lambda_{i}-1+N-\mu_{j}}{\lambda_{i}-1}
\right)_{1\leq i,j\leq n}.
\end{equation*}
Letting $N\rightarrow\infty$,
we obtain the desired identity \eqref{eq:STGF}.
\end{demo}
%
%
%
%
%% ------------------------- %%
%% Proof of Theorem
%% ------------------------- %%
%
%
\begin{demo}{Proof of Theorem~\ref{th:ShiftedTableaux}}
In fact \eqref{eq:ShiftedRPP-GF} is equivalent to the following identity from \eqref{eq:GF-rel}:
\begin{align}
&\sum_{\mu\in\mathcal{P}_{n}'}
\omega_{r-2(m-2n)-1}'(\mu)\,
\GF{\mathcal{T}_{\lambda,\mu}}
\nonumber\\&=
a^{n(n-1)}
q^{mn+n(n-2)(4n-1)/3+n(n-1)r}
\left\{
\frac{(abq^2;q)_{\infty}}{(aq;q)_{\infty}}
\right\}^{n}
\nonumber\\&\times
\prod_{k=1}^{2n}
\frac{1}{(q;q)_{k+m-2n-1}}
\prod_{k=1}^{n-1}(bq;q)_{2k}
\prod_{k=1}^n\frac{(q;q)_{2k-1}(aq;q)_{2k+r-1}}
{(abq^2;q)_{2(k+n)+r-3}}.
\label{eq:ShiftedTableaux-GF}
\end{align}
Now the entries of the skew-symmetric matrix in Corollary~\ref{conj:01}
can be written as
\begin{align*}
a^{i}_{j}
%=(q^{i-1}-q^{j-1})\frac{(aq;q)_{i+j+r-2}}{(abq^2;q)_{i+j+r-2}}
=(q^{i-1}-q^{j-1})\frac{(aq;q)_{\infty}}{(abq^2;q)_{\infty}}
\cdot
\frac{(abq^{i+j+r};q)_{\infty}}{(aq^{i+j+r-1};q)_{\infty}}.
\end{align*}
If we apply the $q$-binomial theorem \cite[(1.3.2)]{GR}
\begin{equation}
\sum_{k=0}^{\infty}\frac{(a;q)_{k}}{(q;q)_{k}}x^k=\frac{(ax;q)_{\infty}}{(x;q)_{\infty}},
\label{eq:q-binomial-theorem}
\end{equation}
then we obtain
\begin{align*}
a^{i}_{j}
&=(q^{i-1}-q^{j-1})\frac{(aq;q)_{\infty}}{(abq^2;q)_{\infty}}
\sum_{k=0}^{\infty}\frac{(bq;q)_{k}}{(q;q)_{k}}\left(aq^{i+j+r-1}\right)^k
\\
&=q^{-s}\frac{(aq;q)_{\infty}}{(abq^2;q)_{\infty}}
\sum_{k=0}^{\infty}\frac{(bq;q)_{k}}{(q;q)_{k}}\left(aq^{r-2s+1}\right)^k
\begin{vmatrix}
q^{(i+s-1)(k+1)}&q^{(i+s-1)k}\\
q^{(j+s-1)(k+1)}&q^{(j+s-1)k}
\end{vmatrix},
%\label{eq:sum}
\end{align*}
where $s$ is any integer.
%We utilize this identity for cetain enumeration in the following subsections.
Hence we obtain
\begin{align}
&
\sum_{k=0}^{\infty}\frac{(bq;q)_{k}}{(q;q)_{k}}\left(aq^{r-2s+1}\right)^k
\begin{vmatrix}
\frac{q^{(i+s-1)k}}{(q;q)_{i+s-2}}&\frac{q^{(i+s-1)(k+1)}}{(q;q)_{i+s-2}}\\
\frac{q^{(j+s-1)k}}{(q;q)_{j+s-2}}&\frac{q^{(j+s-1)(k+1)}}{(q;q)_{j+s-2}}
\end{vmatrix}
\nonumber\\&=
-q^{s}\frac{(abq^2;q)_{\infty}}{(aq;q)_{\infty}}\cdot
\frac{q^{i-1}-q^{j-1}}{(q;q)_{i+s-2}(q;q)_{j+s-2}}\cdot
\frac{(aq;q)_{i+j+r-2}}{(abq^2;q)_{i+j+r-2}}.
\label{eq:sum}
\end{align}
Set $t^{i}_{j}$ and $\alpha_{j}$ to be
\[
t^{i}_{j}=\frac{q^{(i+s-1)j}}{(q;q)_{i+s-2}},
\qquad
\alpha_{j}=\frac{(bq;q)_{j}}{(q;q)_{j}}\left(aq^{r-2s+1}\right)^{j}
\]
for $i\geq1$ and $j\geq0$.
Let $B=(\beta^{i}_{j})_{i,j\geq0}$ be the skew-symmetric matrix defined by
\[
\beta^{i}_{j}=\begin{cases}
\alpha_{i}
&\text{ if $j=i+1$ for $i=0,1,\dots$,}\\
-\alpha_{j}
&\text{ if $i=j+1$ for $j=0,1,\dots$,}\\
0&\text{ otherwise.}
\end{cases}
\]
If we take $T=(t^{i}_{j})_{1\leq i\leq 2n,\,0\leq j}$ 
and $B=(\beta^{i}_{j})_{i,j\geq0}$ 
in Theorem~\ref{th:msf},
then,
by Proposition~\ref{prop:subpfaffian},
we obtain 
\begin{align}
&\sum_{\mu}
\left(aq^{r-2s+1}\right)^{\sum_{k=1}^{n}\mu_{2k-1}}
\prod_{k=1}^{n}\frac{(bq;q)_{\mu_{2k-1}}}{(q;q)_{\mu_{2k-1}}}
\cdot
\det T^{[2n]}_{\mu}
\nonumber\\&
=(-1)^{n}q^{ns}\left\{
\frac{(abq^2;q)_{\infty}}{(aq;q)_{\infty}}
\right\}^{n}
\prod_{k=1}^{2n}\frac1{(q;q)_{k+s-2}}
\cdot
\Pf\left(a^{i}_{j}\right)_{1\leq i,j\leq 2n},
\label{eq:binom-expansion}
\end{align}
where the sum on the left-hand side runs over 
all $2n$-tuples $\mu=(\mu_{1},\mu_{2},\dots,\mu_{2n})$ of integers
such that $0\leq\mu_{1}<\mu_{2}<\dots<\mu_{2n}$ and $\mu_{2}=\mu_{1}+1$,$\dots$,
$\mu_{2n}=\mu_{2n-1}+1$.
Now we take the shape $\lambda=(m,m-1,\dots,m-2n+1)$ for positive integers $m,n$
such that $m\geq2n$ in Lemma~\ref{lem:STGF}.
Then \eqref{eq:STGF} implies that
$\GF{\mathcal{T}_{\lambda,\mu}}$ equals
\begin{equation}
\det\left(
\frac{q^{(m-i+1)\mu_{j}}}{(q;q)_{m-i}}
\right)_{1\leq i,j\leq 2n}
=(-1)^{n}
\det\left(
\frac{q^{(i+m-2n)\mu_{j}}}{(q;q)_{i+m-2n-1}}
\right)_{1\leq i,j\leq 2n}
\label{eq:GF-ST}
\end{equation}
for a profile $\mu=(\mu_1,\dots,\mu_{2n})$
where $0\leq\mu_{1}<\dots<\mu_{2n}$.
Here the right-hand side is obtained from the left-hand side by reversing
the order of row indices.
Let $s=m-2n+1$.
If we substitute \eqref{eq:GF-ST} into \eqref{eq:binom-expansion}
and use Corollary~\ref{conj:01},
then we obtain
\begin{align*}
&\sum_{\mu}
\left(aq^{r-2s+1}\right)^{\sum_{k=1}^{n}\mu_{2k-1}}
\prod_{k=1}^{n}\frac{(bq;q)_{\mu_{2k-1}}}{(q;q)_{\mu_{2k-1}}}
\cdot
\GF{\mathcal{T}_{\lambda,\mu}}
\\&
=q^{ns}\left\{
\frac{(abq^2;q)_{\infty}}{(aq;q)_{\infty}}
\right\}^{n}
\prod_{k=1}^{2n}\frac1{(q;q)_{k+s-2}}
\\&\times
a^{n(n-1)}q^{n(n-1)(4n+1)/3+n(n-1)r}\prod_{k=1}^{n-1}(bq;q)_{2k}
\prod_{k=1}^n\frac{(q;q)_{2k-1}(aq;q)_{2k+r-1}}
{(abq^2;q)_{2(k+n)+r-3}}.
\end{align*}
This proves \eqref{eq:ShiftedTableaux-GF}.
If we put $\mu=\nu+\epsilon_{2n}$ and 
use the fact that $n(\lambda)=mn(2n-1)-\frac13n(2n-1)(4n-1)$,
then we can prove \eqref{eq:ShiftedRPP-GF} by a direct computation.
The other identities can be proven similarly.
The details are left to the reader.
\end{demo}
Theorem~\ref{th:ShiftedTableaux} treats only reverse plane partitions
whose number of rows is even.
We obtain the case where the number of rows equals $2n-1$
from Corollary~\ref{cor:byproduct}.
Let $\check A=({\check a}^{i}_{j})_{0\leq i,j}$ be the skew-symmetric matrix whose $(i,j)$-entry
for $0\leq i<j$ equals
\[
{\check  a}^{i}_{j}
=\begin{cases}
\frac{(aq;q)_{j+r-1}}{(abq^2;q)_{j+r-1}}
&\text{ if $i=0$ and $j\geq1$,}\\
(q^{i-1}-q^{j-1})\frac{(aq;q)_{i+j+r-2}}{(abq^2,q)_{i+j+r-2}}
&\text{ if $1\leq i<j$.}
\end{cases}
\]
Then it is easy to see that Corollary~\ref{cor:byproduct}
implies that for $n\geq1$
\begin{align}
&\Pf\left({\check a}^{i}_{j}\right)_{0\leq i,j\leq 2n-1}
=a^{(n-1)^2}q^{n(n-1)(4n-5)/3+(n-1)^2r}
\nonumber\\&\qquad\times
\frac{(aq;q)_{r}}{(abq^2;q)_{r}}
\prod_{k=1}^{n-1}
\frac
{(q;q)_{2k}(aq;q)_{2k+r}(bq;q)_{2k-1}}
{(abq^2;q)_{4k+r-1}(abq^{2k+r};q)_{2k-1}}.
\label{eq:pf-byproduct-b}
\end{align}
Let ${\check P}_{n,r}(a,b;q)$ denote the right-hand side of \eqref{eq:pf-byproduct-b}.
Then, more generally, from \eqref{eq:pf-general3} and \eqref{eq:pf-general4} we derive
for $n\geq2$
\begin{align}
&\Pf\left({\check A}_{[0,2n-2],m-1}\right)
=
%\nonumber\\&\times
\frac{(q^{m-2n+1};q)_{2n-2}(aq^{2n+r-1};q)_{m-2n}}
{(q;q)_{2n-2}(abq^{4n+r-3};q)_{m-2n}}
{\check P}_{n,r}(a,b;q),
\label{eq:pf-general-b3}
\\
&\Pf\left({\check A}_{[0,2n-3],2n-1,m-1}\right)
=
q\cdot\frac{
(q^{m-2n};q)_{1}
(q^{m-2n+2};q)_{2n-3}
(aq^{2n+r-1};q)_{m-2n}
}
{
(q;q)_{2n-2}
(abq^{4n+r-5};q)_{1}
(abq^{4n+r-3};q)_{m-2n+1}
}
\nonumber\\&\qquad\qquad\qquad\qquad\qquad\times
f(2n-1,m-1,r){\check P}_{n,r}(a,b;q)
.
\label{eq:pf-general-b4}
\end{align}
As an application of \eqref{eq:pf-byproduct-b}, \eqref{eq:pf-general-b3} and \eqref{eq:pf-general-b4},
we can derive a similar identities in the case where the number of rows of the shapes is odd.
Before we state our theorem we need a few definitions.
%
%
%% ------------------------------------- %%
%% Profile and weight
%% ------------------------------------- %%
Fix positive integers $n$ and $t$ such that $1\leq t\leq n$.
Let $\mathcal{Q}^{(t)}_{n}$ denote the set of profiles
$\nu=(\nu_{1},\dots,\nu_{2n-1})$ such that $0\leq\nu_{1}\leq\nu_{2}\leq\dots\leq\nu_{2n-1}$,
$\nu_{2k}=\nu_{2k-1}$ for $k=1,\dots,t-1$
and
$\nu_{2k+1}=\nu_{2k}$ for $k=t,\dots,n-1$.
For $\nu\in\mathcal{Q}^{(t)}_{n}$ and $x,y\in\mathbb{Z}$ we let
\begin{align}
\psi_{x,y}^{(t)}(\nu)
&=\left(aq^{x}\right)^{(|\nu|-\nu_{2t-1})/2}\left(aq^{y}\right)^{\nu_{2t-1}}
\nonumber\\&\times
\frac{(bq^{2t-1};q)_{\nu_{2t-1}}}{(q^{2t-1};q)_{\nu_{2t-1}}}
\prod_{k=1}^{t-1}\frac{(bq^{2k};q)_{\nu_{2k-1}-1}}{(q^{2k};q)_{\nu_{2k-1}-1}}
\prod_{k=t}^{n-1}\frac{(bq^{2k};q)_{\nu_{2k}}}{(q^{2k};q)_{\nu_{2k}}},
\label{eq:profile-weight2}
\end{align}
where $|\nu|=\sum_{k=1}^{2n-1}\nu_{k}$.
Then we obtain the following theorem 
from \eqref{eq:pf-byproduct-b}, \eqref{eq:pf-general-b3} and \eqref{eq:pf-general-b4}.
%
%
%% -------------------------------------------------- %%
%% Theorem (Weighted enumration of shifted tableaux)
%% -------------------------------------------------- %%
%
%
\begin{theorem}
\label{th:ShiftedTableauxOdd}
Let $r$ be an integer.
For any positive integers $m$ and $n$ such that $m\geq2n-1$,
we fix the shape $\lambda=(m,m-1,\dots,m-2n+2)$ of length $2n-1$.
Then we have
\begin{align}
&\sum_{t=1}^{n}
\left(aq^{r+1}\right)^{t-1}
\frac{(bq;q)_{2(t-1)}}{(q;q)_{2(t-1)}}
\sum_{\nu\in\mathcal{Q}^{(t)}_{n}}
\psi_{r-2(m-2n)-3,r-(m-2n)-1}^{(t)}(\nu)\,
\GF{\mathcal{R}_{\lambda,\nu}}
\nonumber\\&=
\prod_{k=1}^{2n-1}
\frac{(q;q)_{k-1}}{(q;q)_{k+m-2n}}
\cdot
R_{n,r}(a,b;q),
\label{eq:ShiftedRPP-GF-odd}
\end{align}
where
\begin{equation*}
R_{n,r}(a,b;q)
=\left\{
\frac{(abq^2;q)_{\infty}}{(aq;q)_{\infty}}
\right\}^{n}
\cdot
\frac{(aq;q)_{r}}{(abq^2;q)_{r}}
\cdot
\prod_{k=1}^{n-1}
\frac{(aq;q)_{2k+r}}
{(abq^2;q)_{4k+r-1}(abq^{2k+r};q)_{2k-1}}.
\end{equation*}
More generally,
if $\lambda=(l,m-1,m-2,m-3,\dots,m-2n+2)$ where $l\geq m$ and $n\geq2$,
then we have
\begin{align}
&\sum_{t=1}^{n}
\left(aq^{r+1}\right)^{t-1}
\frac{(bq;q)_{2(t-1)}}{(q;q)_{2(t-1)}}
\sum_{\nu\in\mathcal{Q}^{(t)}_{n}}
\psi_{r-2(m-2n)-3,r-(m-2n)-1}^{(t)}(\nu)\,
\GF{\mathcal{R}_{\lambda,\nu}}
\nonumber\\&=
\frac{
\prod_{k=1}^{2n-1}(q;q)_{k-1}
}{
(q;q)_{l-1}\prod_{k=1}^{2n-2}(q;q)_{k+m-2n}
}
%%\nonumber\\&\times
\cdot
\frac{
(q^{l-m+1};q)_{2n-2}(aq^{2n+r-1};q)_{l-m}
}{
(q;q)_{2n-2}(abq^{4n+r-3};q)_{l-m}
}
R_{n,r}(a,b;q),
\label{eq:ShiftedRPP-GF-odd2}
\end{align}
and
if $\lambda=(l,m,m-2,m-3,\dots,m-2n+2)$ where $l>m$ and $n\geq2$,
 then we have
\begin{align}
&\sum_{t=1}^{n}
\left(aq^{r+1}\right)^{t-1}
\frac{(bq;q)_{2(t-1)}}{(q;q)_{2(t-1)}}
\sum_{\nu\in\mathcal{Q}^{(t)}_{n}}
\psi_{r-2(m-2n)-3,r-(m-2n)-1}^{(t)}(\nu)\,
\GF{\mathcal{R}_{\lambda,\nu}}
\nonumber\\&=
\frac{
f(2n-1,l-m+2n-1,r)
\prod_{k=1}^{2n-1}(q;q)_{k-1}
}{
(q;q)_{l-1}(q;q)_{m-1}\prod_{k=1}^{2n-3}(q;q)_{k+m-2n}
}
\nonumber\\&\qquad\times
\frac{
(q^{l-m};q)_{1}
(q^{l-m+2};q)_{2n-3}
(aq^{2n+r-1};q)_{l-m}
}
{
(q;q)_{2n-2}
(abq^{4n+r-5};q)_{1}
(abq^{4n+r-3};q)_{l-m+1}
}
R_{n,r}(a,b;q).
\label{eq:ShiftedRPP-GF-odd3}
\end{align}
\end{theorem}
To prove this theorem,
define a matrix $T=(t^{i}_{j})_{i\geq0,\,j\geq-1}$ by
\[
t^{i}_{j}=\begin{cases}
1&\text{ if $i=0$ and $j=-1$,}\\
\frac{q^{(i+s-1)j}}{(q;q)_{i+s-2}}
&\text{ if $i\geq1$ and $j\geq0$,}\\
0&\text{ otherwise,}
\end{cases}
\]
and a skew-symmetric matrix $B=(\beta^{i}_{j})_{-1\leq i<j}$ by
\[
\beta^i_j=
\begin{cases}
\frac{(bq;q)_j}{(q;q)_j}(aq^{r-s+1})^j &
\mbox{ if $i=-1$ and $j\geq 0$,}\cr
\frac{(bq;q)_i}{(q;q)_i}(aq^{r-2s+1})^i &
\mbox{ if $0\leq i<j$ and $i+1=j$,}\cr
0 &\mbox{otherwise,}
\end{cases}
\]
where $s=m-2n+2$.
A similar reasoning as in the proof of Theorem~\ref{th:ShiftedTableaux}
works to prove these identities.
We omit the details.
%
%
%
%
%

%% file: tagawa06.tex
%
%

%%%% ---------------< ooo >--------------- %%%%
%%%%
%%%% Section 6
%%%%
%%%% ---------------< ooo >--------------- %%%%

\section{Open problems}
\label{sec:open-problems}

In this section we formulate several conjectures for the Pfaffians of  certain sequences related to Catalan numbers based on the computer experiments. 
%% ------------ < *** > ------------ %%
%% Al-Salam-Carlitz
%% ------------ < *** > ------------ %%
The \defterm{Al-Salam-Carlitz polynomials} \cite{AC} are defined by
\[
U_{n}^{(a)}(x;q)
=(-1)^{n}q^{\binom{n}2}{}_{2}\phi_{1}
\biggl(
{{q^{-n},ax^{-1}}\atop{0}}\,;\,q,qx
\biggr).
\]
Let $L$ be the linear functional with respect to which
$U_{n}^{(a)}(x;q)$ are orthogonal.
Then the $n$th moment has the expression \cite{GR,Ki,KLS}:
\[
G_{n}(a;q)=L\left(x^{n}\right)=\sum_{k=0}^{n}\left[{{n}\atop{k}}\right]_{q}a^{k},
\]
where $\left[{{n}\atop{k}}\right]_{q}=\frac{(q;q)_{n}}{(q;q)_{k}(q;q)_{n-k}}$.
%
%
%% ------------ < *** > ------------ %%
%% Conjecture
%% ------------ < *** > ------------ %%
%The following is an independent conjecture related to another monent sequence.
\begin{conjecture}
Let $n\geq1$ be an integer.
Then the following identities would hold:
\begin{align}
&\Pf\biggl((q^{i-1}-q^{j-1})G_{i+j-3}(a;q)\biggr)_{1\leq i,j\leq 2n}
\nonumber\\&=
a^{n(n-1)}q^{\frac13\lfloor n/2\rfloor(16\lfloor n/2\rfloor^2-1)-(-1)^n4\lfloor n/2\rfloor^2-2\lfloor n/2\rfloor\cdot\lfloor (n-1)/2\rfloor}\prod_{k=1}^{n}(q;q)_{2k-1},
\\
&\Pf\biggl((q^{i-1}-q^{j-1})G_{i+j-2}(a;q)\biggr)_{1\leq i,j\leq 2n}
\nonumber\\&=
a^{n(n-1)}q^{\frac13\lfloor n/2\rfloor(16\lfloor n/2\rfloor^2-1)-(-1)^n4\lfloor n/2\rfloor^2}\prod_{k=1}^{n}(q;q)_{2k-1}
\sum_{k=0}^{n}
q^{\lfloor (n-2k)^2/2\rfloor}\left[{{n}\atop{k}}\right]_{q^2}a^{k}.
\end{align}
Here $\lfloor x\rfloor$ denotes the largest integer which is not greater than $x$,
and we use the convention that $G_{-1}(a;q)=0$ which can in fact be assigned to any value.
\end{conjecture}
%
%
%% ------------ < *** > ------------ %%
%% Motzkin, Delannoy, Schr\"oder, Narayana
%% ------------ < *** > ------------ %%
There are several well-known numbers related to lattice path enumeration (\cite{A,Sta}).
Let $M_{n}=\sum_{k=0}^{n}\binom{n}{2k}C_{k}$ denote the \defterm{Motzkin numbers},
 $D_{n}=\sum_{k=0}^{n}\binom{n}{k}\binom{n+k}{k}$ the \defterm{central Delannoy numbers},
and $S_{n}=\sum_{k=0}^{n}\binom{n+k}{2k}C_{k}$ \defterm{Schr\"oder numbers}.
%%We use the convention that $M_{-1}=D_{-1}=0$.
%%(In fact we can assign  any value.)
%%The  $S_n$ are given by the recurrence relation
%%$S_n=S_{n-1}+\sum_{k=0}^{n-1}S_kS_{n-1-k}$, 	
%%where $S_{0}=1$.
%
Finally, the number 
 $N(n, k)=\frac1{n}\binom{n}{k}\binom{n}{k-1}$ 
%% ($n = 1, 2, 3\dots$, $1\leq k\leq n$)
is known as a  \defterm{Narayana number}, and 
\begin{align*}
N_{n}(a)=\sum_{k=0}^{n}\frac1{n}\binom{n}{k}\binom{n}{k-1}a^k
\end{align*}
is known as the $n$th \defterm{Narayana polynomial},
which is the moment sequence of a {\sl generalized Chebyshev polynomials} of the first kind.
Here we use the convention that $N_{0}(a)=1$.
%
%
%
%% ---------------------------------------- %%
%% Conjecture
%% ---------------------------------------- %%
\begin{conjecture}
Let $n\geq1$ be an integer.
Then the following identities would hold:
\begin{align}
&\Pf\biggl((j-i)M_{i+j-3}\biggr)_{1\leq i,j\leq 2n}
=\prod_{k=0}^{n-1}(4k+1),
\label{eq:conj-M_n}
\\
&\Pf\biggl((j-i)D_{i+j-3}\biggr)_{1\leq i,j\leq 2n}
=2^{n^2-1}(2n-1)\prod_{k=1}^{n-1}(4k-1),
\label{eq:conj-D_n}
\\
&\Pf\biggl((j-i)S_{i+j-2}\biggr)_{1\leq i,j\leq 2n}
=2^{n^2}\prod_{k=0}^{n-1}(4k+1),
\label{eq:conj-S_n}
\\
&\Pf\biggl((j-i)N_{i+j-2}(a)\biggr)_{1\leq i,j\leq 2n}
=a^{n^2}\prod_{k=0}^{n-1}(4k+1).
\label{eq:conj-N_n}
\end{align}
\end{conjecture}
Note that $C_{n}=N_{n}(1)$, $S_{n}=N_{n}(2)$
and 
\[
M_{n-1}=\left(\frac{1-\sqrt{-3}}{2}\right)^{n+1}
N_n\left(\frac{-1+\sqrt{-3}}{2}\right).
\]
Hence, if one could prove \eqref{eq:conj-N_n},
then one would have proven  \eqref{eq:conj-M_n} and \eqref{eq:conj-S_n} as corollaries.

\smallskip
Let $a_{n}=\frac1{2n+1}\binom{3n}{n}=\frac1{3n+1}\binom{3n+1}{n}$.
In \cite{GX} Gessel and Xin prove that $\det(a_{i+j-1})_{1\leq i,j\leq n}$ equals
the number of $(2n+1)\times(2n+1)$ alternating sign matrices
that are invariant under vertical reflection.
We propose the following conjecture concerning this sequence.
\begin{conjecture}
Let $a_{n}$ be as above.
Then the following identity would hold:
\begin{equation}
\Pf\left((j-i)a_{i+j-1}\right)_{1\leq i,j\leq2n}
=\frac1{2^n}\prod_{k=1}^{n}\frac
{(12k-6)!(4k-3)!(3k-1)!}{(8k-6)!(8k-3)!(3k-2)!}.
\end{equation}
\end{conjecture}

%% file: tagawa07.tex
%
%

%%%% ---------------< ooo >--------------- %%%%
%%%%
%%%% Appendix
%%%%
%%%% ---------------< ooo >--------------- %%%%

%\appendix
\section*{Appendix: Creative telescoping}

In this appendix we state an alternative proof of \eqref{eq:k:odd} and \eqref{eq:k:even} by
 Zeilberger's creative telescoping
\cite{Ko,PWZ}.
In this case the certificates are extremely simple,
and we can check the computation by hand.
We note that one can prove \eqref{eq:add} similarly,
but the certificate for \eqref{eq:add} is a little more complicated.
\par
\medskip
%% ------------ < *** > ------------ %%
%% Zeilberger's creative telescoping
%% ------------ < *** > ------------ %%
By replacing $q^i$ by $c$,
the equations \eqref{eq:k:odd} and \eqref{eq:k:even} are generalized as
\begin{align}
&\sum
a^{k-1}q^{k(k-1)+1}
\cdot
\frac{
(abq^{2k};q)_{1}
(abq^2;q)_{k-2}
(bq,cq^{-k+1},q^{j-k+1};q)_{k-1}
}{
(q,aq,abcq,abq^{j+1};q)_{k}
}
\nonumber\\&\qquad\times
h_{k}(j;a,b,c,q)
%%%\nonumber\\&
=\frac{
(ac,abq^2;q)_{j-1}
}{
(aq,abcq\,;q)_{j-1}
},
\label{eq:k:odd-even}
\end{align}
where
the sum on the left-hand side runs over odd positive integers
or even nonnegative integers,
and $h_{k}(j;a,b,c,q)$ is set to be
\begin{align}
&h_{k}(j;a,b,c,q)=(1-q^{k})(1-aq^{k})
\nonumber\\&\qquad\times
\Bigl\{
q^{-k}(1+abq^{2k})(1+abcq^{j-1})
-ab(1+q)\left(cq^{-1}+q^{j-1}\right)\Bigr\}
\nonumber\\&\qquad
+aq^{k-1}(1-b)
\left(1-cq^{-k}\right)\left(1-q^{j-k}\right)
\left(1-abq^{2k+1}\right).
\label{eq:h_{k}(j;a,b,c,q)}
\end{align}
%%%%%%%%%%%%%%%%
%%%%%%%%%%%%%%%%%%%%%%%%%%%%%%%%%%%%%
%
%
%% 25/11/2010 %% Hence $h_{k}(j;a,b,c,q)$ is a Laurent polynomial of $q^k$ of ldegree $-1$ and degree $3$.
%
%
%
%
%
Let
%% 16/11/2010 %% \begin{align}
%% 16/11/2010 %% F(j,k)&=
%% 16/11/2010 %% a^{k-1}q^{k(k-1)+1}
%% 16/11/2010 %% \cdot
%% 16/11/2010 %% \frac{
%% 16/11/2010 %% (abq^{2k};q)_{1}
%% 16/11/2010 %% (abq^2;q)_{k-2}
%% 16/11/2010 %% (bq,cq^{-k+1},q^{j-k+1};q)_{k-1}
%% 16/11/2010 %% }{
%% 16/11/2010 %% (q,aq,abcq,abq^{j+1};q)_{k}
%% 16/11/2010 %% }
%% 16/11/2010 %% \nonumber\\&\qquad\times
%% 16/11/2010 %% h_{k}(j;a,b,c,q)
%% 16/11/2010 %% \label{eq:F(j,k)}
%% 16/11/2010 %% \end{align}
\begin{align}
F(j,k)&=
a^{k-1}c^{k-1}q^{j(k-1)+1}
\cdot
\frac{
(abq^{2k};q)_{1}
(abq^2;q)_{k-2}
(bq,q/c,q^{1-j};q)_{k-1}
}{
(q,aq,abcq,abq^{j+1};q)_{k}
}
\nonumber\\&\qquad\times
h_{k}(j;a,b,c,q).
\label{eq:F(j,k)}
\end{align}
%
%
%
%% 24/11/2010 %% For a fixed nonnegative integer $N$
%% 24/11/2010 %% and fixed coefficients $c_{n}(j)$,
%% 24/11/2010 %% we let
%% 24/11/2010 %% \begin{align}
%% 24/11/2010 %% T(j,k)=\sum_{n=0}^{N}c_{n}(j)F(j+n,k),
%% 24/11/2010 %% \end{align}
%% 24/11/2010 %% Hence we have
%% 24/11/2010 %% \begin{align*}
%% 24/11/2010 %% &T(j,k)=a^{k-1}c^{k-1}q^{j(k-1)+1}
%% 24/11/2010 %% \cdot
%% 24/11/2010 %% \frac{
%% 24/11/2010 %% (1-abq^{2k})
%% 24/11/2010 %% (abq^2;q)_{k-2}
%% 24/11/2010 %% (bq,q/c,q^{1-j};q)_{k-1}
%% 24/11/2010 %% }{
%% 24/11/2010 %% (q,aq,abcq;q)_{k}
%% 24/11/2010 %% (abq^{j+1};q)_{k+N}
%% 24/11/2010 %% (q^{k-j-N};q)_{N}
%% 24/11/2010 %% }
%% 24/11/2010 %% \\&\times
%% 24/11/2010 %% \sum_{n=0}^{N}
%% 24/11/2010 %% c_{n}(j)q^{n(k-1)}
%% 24/11/2010 %% (abq^{j+1},q^{1-j-n};q)_{n}
%% 24/11/2010 %% (abq^{j+k+n+1},q^{k-j-n};q)_{N-n}
%% 24/11/2010 %% h_{k}(j+n;a,b,c,q).
%% 24/11/2010 %% \end{align*}
Hereafter we use the notation that 
$F^{(o)}(j,k)=F(j,2k-1)$
and
$F^{(e)}(j,k)=F(j,2k-2)$.
Further
we set
$T^{(o)}(j,k)=T(j,2k-1)$
and
$T^{(e)}(j,k)=T(j,2k-2)$,
where
\begin{equation}
T(j,k)=F(j,k)
-\frac { \left( 1-aq^{j} \right)  \left( 1-abcq^j \right)}{ \left( 1-abq^{j+1} \right)  \left( 1-acq^{j-1} \right) }
F(j+1,k).
\end{equation}
Let us define
$P^{(x)}(j,k)$, $Q^{(x)}(j,k)$ and $R^{(x)}(j,k)$ for
$x=o,e$ by
$P^{(o)}(j,k)=P(j,2k-1)$, $Q^{(o)}(j,k)=Q(j,2k-1)$,
$R^{(o)}(j,k)=R(j,2k-1)$,
$P^{(e)}(j,k)=P(j,2k-2)$, $Q^{(e)}(j,k)=Q(j,2k-2)$,
 and $R^{(e)}(j,k)=R(j,2k-2)$,
where
\begin{align*}
P(j,k)
&=a^2c^2q^{2j}(1-abq^{k})
(1-abq^{k+1})
(1-bq^{k})
(1-bq^{k+1})
\\&\qquad\times
(1-q^{k}/c)
(1-q^{k+1}/c)
(1-q^{k-j-1})
(1-q^{k-j}),
\\
%% 24/11/2010 %% \widetilde Q(j,k)
%% 24/11/2010 %% &=(1-q^{k+1})(1-q^{k+2})(1-aq^{k+1})(1-aq^{k+2})(1-abcq^{k+1})
%% 24/11/2010 %% \\&\qquad\times
%% 24/11/2010 %% (1-abcq^{k+2})(1-abq^{j+k+N+1})(1-abq^{j+k+N+2}),
%% 24/11/2010 %% \\
Q(j,k)
&=(1-q^{k+1})(1-q^{k+2})(1-aq^{k+1})(1-aq^{k+2})(1-abcq^{k+1})
\\&\qquad\times
(1-abcq^{k+2})(1-abq^{j+k+2})(1-abq^{j+k+3}),
\\
%% 24/11/2010 %% R(j,k)&=(1-abq^{2k})\,
%% 24/11/2010 %% \sum_{n=0}^{N}
%% 24/11/2010 %% c_{n}(j)(abq^{j+1},q^{1-j-n};q)_{n}
%% 24/11/2010 %% \\&\qquad\qquad\qquad\qquad\times
%% 24/11/2010 %% (abq^{j+k+n+1},q^{k-j-n};q)_{N-n}
%% 24/11/2010 %% h_{k}(j+n;a,b,c,q),
%% 24/11/2010 %% \\
R(j,k)&=(1-abq^{2k})\,
\biggl\{
(abq^{j+k+1},q^{k-j-1};q)_{1}
h_{k}(j;a,b,c,q)
\\&
-\frac { q^{k-1}\left( 1-aq^{j} \right)  \left( 1-abcq^j \right)}{ \left( 1-abq^{j+1} \right)  \left( 1-acq^{j-1} \right) }
(abq^{j+1},q^{-j};q)_{1}
h_{k}(j+1;a,b,c,q)\biggr\}.
\end{align*}
By direct computation, we see that
\begin{align}
\frac{T^{(x)}(j,k+1)}{T^{(x)}(j,k)}
=\frac{
P^{(x)}(j,k)
}{Q^{(x)}(j,k)
}\cdot\frac{R^{(x)}(j,k+1)}{R^{(x)}(j,k)}
\label{eq:ratio}
\end{align}
holds for $x=o,e$.
%% 24/11/2010 %% Then the polynomials $P_{x}(j,k)$ and $Q_{x}(j,k)$ satisfy
%% 24/11/2010 %% \begin{equation}
%% 24/11/2010 %% \gcd\left(P_{x}(j,k),Q_{x}(j,k+l)\right)=1
%% 24/11/2010 %% \end{equation}
%% 24/11/2010 %% for all nonnegative integers $l$.
%
%
%
%% 24/11/2010 %% Hence $P_{x}(j,k)$ and $Q_{x}(j,k)$ are polynomials of $q^k$,
%% 24/11/2010 %% and $\deg P_{x}(j,k)=\deg Q_{x}(j,k)=4$ (resp. $16$)
%% 24/11/2010 %% for $x=\emptyset$ (resp. $x=o,e$).
%% 24/11/2010 %% Meanwhile $R_{x}(j,k)$ is a Laurent polynomial of $q^k$
%% 24/11/2010 %% of ldegree $-1$ (resp. $-2$) and degree $2N+5$ (resp. $4N+10$)
%% 24/11/2010 %% for $x=\emptyset$ (resp. $x=o,e$).
%
%
%
We define $\Lambda^{(x)}(j,k)$ by
\begin{equation}
\Lambda^{(x)}(j,k)=\frac{Q^{(x)}(j,k-1)T^{(x)}(j,k)}{R^{(x)}(j,k)}\,X^{(x)}(j,k)
\label{eq:X(j,k)}
\end{equation}
for $x=o,e$, where
%% 25/11/2010 %% \begin{align}
%% 25/11/2010 %% X(j,k)&=
%% 25/11/2010 %% -\frac {1+ac{q}^{j-1}}{ 1-ac{q}^{j+1} }q^{-k}
%% 25/11/2010 %% \nonumber\\&
%% 25/11/2010 %% +\frac{\left(1+q^2\right)
%% 25/11/2010 %% \left\{q+aq+ac+abcq+{q}^{j+1} \left( a+abc+{a}^{2}bc+abq \right)\right\}
%% 25/11/2010 %% }{
%% 25/11/2010 %% {q}^{2} \left( 1+q \right)  \left( 1-ac{q}^{j+1} \right) 
%% 25/11/2010 %% }
%% 25/11/2010 %% \nonumber\\&
%% 25/11/2010 %% -\frac{a\left(1+q^2\right)
%% 25/11/2010 %% \left\{1+bc+abc+b \left( 1+a+abc \right) {q}^{j+1}\right\}
%% 25/11/2010 %% }{
%% 25/11/2010 %% {q}^{2} \left(1- ac{q}^{j+1} \right) 
%% 25/11/2010 %% }q^k
%% 25/11/2010 %% \nonumber\\&
%% 25/11/2010 %% +\frac{ab\left(1+q^2\right)
%% 25/11/2010 %% \left\{q+aq+ac+abcq+{q}^{j+1} \left( a+abc+{a}^{2}bc+abq \right)\right\}
%% 25/11/2010 %% }{
%% 25/11/2010 %% {q}^{2} \left( 1+q \right)  \left( 1-ac{q}^{j+1} \right) 
%% 25/11/2010 %% }q^{2k}
%% 25/11/2010 %% \nonumber\\&
%% 25/11/2010 %% -\frac {{a}^{2}{b}^{2} \left( 1+ac{q}^{j-1} \right) }{1-ac{q}^{j+1} }
%% 25/11/2010 %% q^{3k},
%% 25/11/2010 %% \end{align}
%% 25/11/2010 %% and 
$X^{(o)}(j,k)=X(j,2k-1)$
and $X^{(e)}(j,k)=X(j,2k-2)$,
with
\begin{equation}
X(j,k)
=-\frac{q^{-k}}{1-acq^{j-1}}.
\end{equation}
\begin{lemma}
Let $\Lambda^{(x)}(j,k)$ be as above for $x=o,e$.
Then we have 
\begin{align}
T^{(x)}(j,k)=\Lambda^{(x)}(j,k+1)-\Lambda^{(x)}(j,k)
\label{eq:recurrence}
\end{align}
for $x=o,e$.
\end{lemma}
\begin{demo}{Proof}
If one use \eqref{eq:ratio} and \eqref{eq:X(j,k)}
then \eqref{eq:recurrence} reduces to the following identity:
\begin{align}
&P^{(x)}(j,k)X^{(x)}(j,k+1)- Q^{(x)}(j,k-1) X^{(x)}(j,k)
= R^{(x)}(j,k),\label{eq:Gosper-(4.5)}
\end{align}
This can be checked by direct computation.
\end{demo}
Because of $\Lambda^{(x)}(j,1)=0$ for $x=o,e$,
by summing \eqref{eq:recurrence} over all positive integers,
we obtain
\begin{align}
\sum_{k\geq1}F^{(x)}(j+1,k)
=
\frac
{ \left( 1-abq^{j+1} \right)  \left( 1-acq^{j-1} \right) }
{ \left( 1-aq^{j} \right)  \left( 1-abcq^j \right)}
\sum_{k\geq1}F^{(x)}(j,k)
\label{eq:rec-j}
\end{align}
for $x=o,e$.
Since $F(1,0)=F(1,1)=1$ and $F(1,k)=0$ for $k\geq2$,
we have
\begin{equation}
\sum_{k\geq1}F^{(x)}(1,k)=1
\end{equation}
for $x=o,e$.
Hence we obtain the desired identity \eqref{eq:k:odd-even} from \eqref{eq:rec-j}.
This gives the second proof of Theorem~\ref{conj:decomposition}
and Theorem~\ref{th:byproduct}.